\documentclass[a4paper,12pt]{article}

\usepackage{amsmath,latexsym,amssymb}
\usepackage[german,english]{babel}
\newtheorem{thm}{Theorem}[section]
\newtheorem{prop}[thm]{Proposition}
\newtheorem{lemma}[thm]{Lemma}

\newtheorem{main}{Main Theorem}

\newcommand{\proof}[1][]{{\it Proof#1: }}
\newcommand{\qed}[1][3mm]{\hspace*{\fill} $\Box$ \vspace{#1}}

\newcommand{\ZZ}{{\mathbf Z}}
\newcommand{\CC}{{\mathbf C}}
\newcommand{\PP}{{\mathbf P}}

\newcommand{\del}{\partial}
\newcommand{\ON}{\operatorname}

\renewcommand{\a}{\alpha}
\renewcommand{\b}{\beta}

\renewcommand{\d}{\delta}

\newcommand{\z}{\zeta}

\newcommand{\varth}{\vartheta}

\renewcommand{\l}{\lambda}

\newcommand{\s}{\sigma}
\newcommand{\oo}{\omega}

\newcommand{\tto}{\longrightarrow}

\newcommand{\surj}{\to\!\!\!\!\!\!\!\!\tto}

\newcommand{\inv}{^{^{-1}}}

\newcommand{\afami}{{\cal A}}
\newcommand{\bfami}{{\cal B}}
\newcommand{\cfami}{{\cal C}}
\newcommand{\dfami}{{\cal D}}

\newcommand{\ffami}{{\cal F}}

\newcommand{\gfami}{{\cal G}}
\newcommand{\hfami}{{\cal H}}

\newcommand{\rfami}{{\cal R}}

\newcommand{\ufami}{{\cal U}}
\newcommand{\vfami}{{\cal V}}

\newcommand{\ibold}{{\boldsymbol{i}}}
\newcommand{\jbold}{{\boldsymbol{j}}}
\newcommand{\kbold}{{\boldsymbol{k}}}
\newcommand{\ubold}{{\boldsymbol{u}}}
\newcommand{\xbold}{{\boldsymbol{x}}}
\newcommand{\nubold}{{\boldsymbol{\nu}}}

\renewcommand{\ubold}{u}
\renewcommand{\xbold}{x}
\renewcommand{\nubold}{\nu}

\newcommand{\cutoff}[1]{}
\newcommand{\labell}[1]{\label{#1}}
\newcommand{\reff}[1]{\ref{#1}}

\newcommand{\br}{\ON{Br}}

\newcommand{\psl}{\ON{PSL_2\ZZ}}

\newcommand{\slz}{\ON{SL_2\ZZ}}

\newcommand{\dip}{p_{n,d}}
\newcommand{\bip}{q_{n,d}}
\newcommand{\pn}{\mathbf{P}^n}
\newcommand{\euhyp}{e_{n,d}}
\newcommand{\euci}{e_{n;d,d}}
\newcommand{\pe}{\mathbf{P}^1}
\newcommand{\pp}{\mathbf{P}}
\newcommand{\lep}{\ell_{n,d}}
\newcommand{\degz}{\deg_z}
\newcommand{\wdeg}{\ON{w-deg}}
\newcommand{\degv}{\deg_v}

\newcommand{\hnd}{\hfami_{n,d}}

\newcommand{\ind}{I_{n,d}}

\newcommand{\psym}{\PP\ON{Sym}^d\CC^{n+1}}
\newcommand{\sym}{\ON{Sym}^d\CC^{n+1}}
\newcommand{\und}{\ufami_{n,d}}
\newcommand{\uell}{\ufami_{2,3}}
\newcommand{\upezz}{\ufami_{3,3}}

\newcommand{\diff}{\ON{Diff}^+}
\newcommand{\diffo}{\ON{Diff}^o}
\newcommand{\aut}{\ON{Aut}}

\newcommand{\ncount}{\kappa}
\newcommand{\vect}{{\sf v}}
\newcommand{\kord}{\prec_\ncount}

\begin{document}

\title{Fundamental Groups of Spaces of Smooth Projective Hypersurfaces
\footnote{Project 13}}

\author{Michael L\"onne}

\maketitle

\begin{abstract}
We investigate the complement of the discriminant in the projective space
$\PP\ON{Sym}^d\CC^{n+1}$ of polynomials defining hypersurfaces of degree $d$
in $\PP^n$.
Following the ideas of Zariski we are able to give a presentation for the fundamental
group of the discriminant complement which generalises the well-known
presentation in case $n=1$, i.e.\ of the spherical braid group on $d$ strands.

In particular our argument proceeds by a geometric analysis of the discriminant
polynomial as proposed in \cite{bessis} and draws on results and methods
from \cite{habil} addressing a comparable problem for any versal unfolding
of Brieskorn Pham singularities.
\end{abstract}

%
%
%
%
%

\section{Introduction}

The primary objects for this paper are hypersurfaces of degree $d$ in
projective space $\pn$ which we equip with homogeneous coordinates
$(x_0:...:x_n)$.
They are assembled in the universal hypersurface $\hnd$ which is a hypersurface
in $\PP^n\times\psym$.

Here $\sym\cong\CC[x_0,...,x_n]_d$ can be understood as the vector space of
all homogeneous polynomials of degree $d$ in $n+1$ variables, which is of dimension
$N=\binom{n+d}{d}$.
On $\CC[x_0,...,x_n]_d$ we introduce coordinates $u_\nubold$ with respect to the monomial
basis $x_0^{\nu_0}\cdots x_n^{\nu_n}$, where the multiindex $\nubold$ is taken from the
multiindex set of non-negative integers $\vfami=\{(\nu_0,...,\nu_n)\,|\,\sum \nu_i=d\}$.

The equation of $\hnd$ then reads in multiindex notation
$$
\sum_{\nubold\in\vfami} u_\nubold \xbold^\nubold \quad=\quad 0.
$$

Since projection of $\hnd$ to the factor $\PP^{N-1}=\psym$ has singular values exactly along the
discriminant
$$
\dfami_{n,d} \quad = \quad \{\ubold\in\PP^{N-1}| \hfami_\ubold \text{ is singular}\}
$$
it is regular over the complement $\und$, where it is thus the projection of a locally trivial fibration
with fibre diffeomorphic to a smooth hypersurface of degree $d$ in $\PP^n$.


The problem we want to address in this paper is to give a finite presentation of the fundamental group
$\pi_1(\und)$ which is geometrically distinguished.
It owes very much to Zariski, who was in fact the first to study it and to settle the case $n=1$.
We cite his result to prepare the stage:

\begin{thm}[Zariski \cite{z} and Fadell, van Buskirk \cite{fvb}]
\hspace*{\fill}\phantom{x}\\
The fundamental group $\pi_1(\ufami_{1,d})$
is finitely presented by generators
$\s_i$, $1\leq i<d$ and
\begin{enumerate}
\item
$\s_i\s_j=\s_j\s_i$, if $|i-j|>1$, $1\leq i,j<d$,
\item
$\s_i\s_{i+1}\s_i=\s_{i+1}\s_i\s_{i+1}$, if $1\leq i<d-1$,
\item
$\s_1\cdots\s_{d-2}\s_{d-1}^2\s_{d-2}\cdots\s_1=1$.
\end{enumerate}
\end{thm}

Our generalisation of his result may be stated as follows.
(For a more precise statement see section \reff{results}.)

\begin{main} 
There is a graph $\Gamma_{n,d}$ with vertex set $V_{n,d}$ and edge set $E_{n,d}\subset{V_{n,d}}^2$
and for $\ncount\in\{0,...,n\}$ there are bijective maps
$$
\ibold_\ncount: \{1,...,(d-1)^n\}\to V_{n,d}
$$
such that $\pi_1(\und)$ is generated by elements $\s_\ibold$ in bijection to the elements $\ibold\in V_{n,d}$
and a complete set of relations is given as follows:
\begin{enumerate}
\item
$\s_\ibold\s_\jbold=\s_\jbold\s_\ibold$ for all $(\ibold,\jbold)\not\in E_{n,d}$,
\item
$\s_\ibold\s_\jbold\s_\ibold=\s_\jbold\s_\ibold\s_\jbold$ for all $(\ibold,\jbold)\in E_{n,d}$,
\item
$\s_\ibold\s_\jbold\s_\kbold\s_\ibold=\s_\jbold\s_\kbold\s_\ibold\s_\jbold$ for all
$(\ibold,\jbold,\kbold)$ such that $(\ibold,\jbold),(\ibold,\kbold),(\jbold,\kbold)\in E_{n,d}$,
\item
for all $\ibold\in V_{n,d}$
$$
\s_\ibold\left(\s_\ibold \prod_{m=1}^{(d-1)^n}\s_{\ibold_0(m)}\right)^{d-1}
\,=\quad \left(\s_\ibold \prod_{m=1}^{(d-1)^n}\s_{\ibold_0(m)}\right)^{d-1} \s_\ibold
$$
\item
$$
\prod_{\ncount=0}^n\:\prod_{m=1}^{(d-1)^n}\s_{\ibold_\ncount(m)}\quad=\quad 1.
$$
\end{enumerate}
\end{main}

Zariski's result is recovered with $\Gamma_{1,d}$ the graph on $d-1$ vertices of type $A_{d-1}$ and the two maps
$\ibold_0$ and $\ibold_1$ enumerating its vertices from right to left and from left to right respectively.
In that case the set of relations $iii)$ is void and the relations in $iv)$ are consequences of the last one.
\\

We want to stress the fact that the relations $i)-iii)$ of our presentation have a distinctive
flavour since they occur in a different setting as well:

If we consider the Fermat polynomial or Brieskorn-Pham polynomial in the variables $x_1,..,x_n$:
$$
x_1^d+\cdots+ x_n^d
$$
we are naturally led to consider a versal unfolding of the isolated hypersurface singularity it defines.
In fact the complement of the discriminant in the unfolding base was shown to have fundamental group
generated as in the theorem but with relations $i)-iii)$ only, \cite{habil}, in terms of a distinguished
Dynkin graph $\Gamma_{n,d}$ associated to the singularity, eg. the following and higher dimensional analoga:
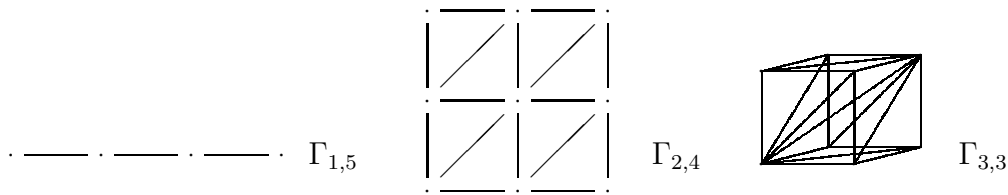
\unitlength=1.2mm

\begin{picture}(130,26)

\put(0,0){\begin{picture}(40,14)(-5,-6)

\put(1.5,0){\line(1,0){7}}
\put(11.5,0){\line(1,0){7}}
\put(21.5,0){\line(1,0){7}}

\put(0,0){\circle*{.3}}
\put(10,0){\circle*{.3}}
\put(20,0){\circle*{.3}}
\put(30,0){\circle*{.3}}

\end{picture}}

\put(38,5){$\Gamma_{1,5}$}

\put(50,0){
\begin{picture}(35,14)(0,-2)

\put(0,0){\circle*{.3}}
\put(0,10){\circle*{.3}}
\put(10,0){\circle*{.3}}
\put(20,0){\circle*{.3}}
\put(20,10){\circle*{.3}}
\put(0,20){\circle*{.3}}
\put(10,20){\circle*{.3}}
\put(20,20){\circle*{.3}}
\put(10,10){\circle*{.3}}

\put(1.5,1.5){\line(1,1){7}}

\put(11.5,1.5){\line(1,1){7}}

\put(1.5,11.5){\line(1,1){7}}

\put(11.5,11.5){\line(1,1){7}}

\put(1.5,0){\line(1,0){7}}
\put(1.5,10){\line(1,0){7}}
\put(11.5,0){\line(1,0){7}}
\put(11.5,10){\line(1,0){7}}
\put(1.5,20){\line(1,0){7}}
\put(11.5,20){\line(1,0){7}}

\put(0,1.5){\line(0,1){7}}
\put(10,1.5){\line(0,1){7}}
\put(20,1.5){\line(0,1){7}}
\put(0,11.5){\line(0,1){7}}
\put(10,11.5){\line(0,1){7}}
\put(20,11.5){\line(0,1){7}}

\end{picture}}

\put(76,5){$\Gamma_{2,4}$}

\put(87,5){\unitlength=.35pt
\begin{picture}(200,200)

\put(0,0){
\begin{picture}(100,100)
\put(0,0){\circle*{4}}
\put(100,0){\circle*{4}}
\put(0,100){\circle*{4}}
\put(100,100){\circle*{4}}

\bezier{110}(0,0)(50,0)(100,0)
\bezier{110}(0,100)(50,100)(100,100)
\bezier{110}(0,0)(0,50)(0,100)
\bezier{110}(100,0)(100,50)(100,100)
\bezier{113}(0,0)(50,50)(100,100)

\bezier{110}(0,0)(36,9)(72,18)
\bezier{110}(100,0)(136,9)(172,18)
\bezier{110}(0,100)(36,109)(72,118)
\bezier{110}(100,100)(136,109)(172,118)

\bezier{113}(0,0)(86,9)(172,18)

\bezier{113}(0,100)(86,109)(172,118)

\bezier{113}(0,0)(36,59)(72,118)

\bezier{113}(100,0)(136,59)(172,118)

\bezier{130}(0,0)(86,59)(172,118)

\end{picture}}

\put(72,18){
\begin{picture}(100,100)
\put(0,0){\circle*{4}}
\put(100,0){\circle*{4}}
\put(0,100){\circle*{4}}
\put(100,100){\circle*{4}}
\bezier{110}(0,0)(50,0)(100,0)
\bezier{110}(0,100)(50,100)(100,100)
\bezier{110}(0,0)(0,50)(0,100)
\bezier{110}(100,0)(100,50)(100,100)
\bezier{113}(0,0)(50,50)(100,100)

\end{picture}}

\end{picture}}

\put(110,5){$\Gamma_{3,3}$}

\end{picture}

We will explain in more detail in section \reff{affine} how this result can be used in the present paper. 
Here we only point out, that the affine singularity also occurs on a hypersurface
which is the projective cone over a smooth hypersurface of one dimension less.

Relation $iv)$ and $v)$ should -- on the contrary -- be regarded as due to degeneration
along the hypersurface $x_0$, complement to the affine open with coordinates $x_1,...,x_n$.
\\


Our interest in the problem stems from several sources:

\paragraph{1.)}
One line of investigating $\pi_1(\und)$ has been via the monodromy map to the
automorphisms of the primitive middle homology of a reference hypersurface,
\cite{bruce}, \cite{b}.
This study was refined by \cite{ct} who instead of the tautological hypersurface
$\hnd$ consider the family of 
associated cyclic branched covers of $\PP^n$ which they showed to capture
much more of $\pi_1(\und)$.

With our presentation at hand it should be worthwhile to investigate the kernel of
each monodromy map.
Based on this study it may be possible to characterise the geometric
features of $\und$ and $\hnd$ which are detected by the various monodromy maps.

It may also be an encouragement to set up and study further
monodromy map -- which may include maps of more geometric than homological
character.
\\

\paragraph{2.)}
In a second approach to $\pi_1(\und)$ the parameter spaces $\uell$ and $\upezz$
of hypersurfaces have been investigated using the wealth of geometric properties
known in case of cubic curves and cubic surfaces, e.g\ detailed descriptions of the
moduli spaces of elliptic curves and of del Pezzo surfaces of degree $3$.
(See \cite{dl} for $\pi_1(\ufami_{2,3})$ and \cite{l,lo,act} for $\pi_1(\ufami_{3,3})$.)

How this analysis is related to our presentation should be the objective of further study.

Looijenga in fact gives alternative presentations not only for $\pi_1(\ufami_{3,3})$ but also for 
$\pi_1(\ufami_{1,d})$ which reflects in some way the toric geometry of $\PP^1$.
Together with results on linear systems in weighted projective spaces we are preparing
it may lay a foundation to proceed to general linear systems on toric varieties.

Reversing the direction of the argument we may try to derive some geometric properties
of moduli spaces related to hypersurfaces from our presentation of fundamental groups.

It seems feasible to expect some results on topological fundamental groups of moduli stacks
which have attracted some interest recently.

One of the outstanding problems should be:

\begin{quote}
To what extend is a classifying map to the discriminant complement determined by the
induced map on fundamental groups?

Is there any rigidity result for algebraic maps or even for analytic maps?
\end{quote}

It is noteworthy in this respect that the moduli space of cubic surfaces has been shown to be an
Eilenberg- MacLane space \cite{act}, so all topological maps to moduli space are determined up to homotopy
by the induced map.

Our results invites an investigation of moduli stacks and their fundamental groups
generalising the fact, that $\slz$ is the fundamental group of the moduli stack
of elliptic curves.

\paragraph{3.)}
A very intriguing observation to our taste is the close relation between notions in algebraic
geometry and their seemingly much more relaxed topological counterparts.

In our setting this becomes manifest when comparing $\ufami_{1,d}$ and the configuration
space $F_d[S^2]$ of $d$-point subsets of the sphere, they are just diffeomorphic and
so are their fundamental groups $\pi_1(\ufami_{1,d})$ and the braid group $\br_d[S^2]$.

It was actually in the setting of configuration spaces that the result of Zariski was reproduced by Fadell/van Buskirk.

Dolgachev and Libgober \cite{dl} proposed to study higher dimensional analogues.
To this end $\und$ should be considered as a `configuration space'
of algebraic submanifolds of $\PP^n$ of fixed topological type.

The appropriate topological space is then different since it should contain all topological
submanifolds isotopic to the algebraic hypersurface.

By results of Cerf \cite{cerf} the space of topological submanifolds of $\PP^n$ isotopic to $H_d$ can be 
identified as a coset space for the group
$\diffo(\PP^n)$ of diffeomorphisms of $\PP^n$ isotopic to the identity
with respect to the subgroup $\diffo(\PP^n,H_d)$ of diffeomorphisms which induce a diffeomorphism
of $H_d$ to itself.

This coset space is the natural topological `configuration space' in higher dimensions
$$
F_{H_d}[\PP^n]\quad=\quad \diffo(\PP^n)/\diffo(\PP^n,H_d)
$$
in analogy to $F_d[S^2]=\diffo(S^2)/\diffo(S^2,\{p_1,...,p_d\})$.

The corresponding quotient map is a fibrations which gives rise to a homotopy exact sequence
$$
\pi_1\diffo(\PP^n)\tto\pi_1(F_{H_d}[\PP^n])\tto\pi_0\diffo(\PP^n,H_d)\tto 1
$$
where of course the middle group should be called the `generalised' braid group of
algebraic hypersurfaces in $\PP^n$

There are various maps interrelating these groups, most notably
\begin{eqnarray}
\pi_1 \aut\PP^n & \tto & \pi_1\diffo\PP^n\\
\pi_1(\und) & \tto & \pi_1(F_{H_d}[\PP^n])\\
\pi_0\diffo(\PP^n,H_d) & \tto & \pi_0\diff(H_d).
\end{eqnarray}
It should be a very rewarding task to deduce some of their properties now that
a presentation for $\pi_1(\und)$ is know.

As a matter of fact we just remind the reader that the last map needs not to be surjective even under the restriction
to orientation preserving maps as Hirose \cite{hirose} shows.

The induced maps
\begin{eqnarray*}
\pi_1(\und) & \tto & \pi_0\diffo(\PP^n,H_d),\\
\pi_1(\und) & \tto & \pi_0\diff(H_d).
\end{eqnarray*}
may be seen as possible geometric monodromy maps alluded to before.

\paragraph{}
Finally we would like to mention the impact of our result in the theory of plane curve complements.
Initialised by Zariski the study of plane curve complements has been revitalised by
Moishezon and many fascinating aspects have been added since, eg. various examples of
homeomorphic curves with non-homeomorphic complements by work of Libgober, Oka, Artal,
Cogolludo among others, and their fundamental role in the study of
projective surfaces and their symplectic analogues, cf. work of Auroux, Donaldson and Katzarkov, Kulikov, Teicher.
\\

Back at its roots Zariski proved that projective hypersurface complements have the same
fundamental group as their restriction to generic planes. But he insists that one should actually use that
result 'the other way' to get fundamental groups for plane curves which are generic sections of some natural
higher dimensional hypersurface. 

From the discriminant $\dfami_{1,d}$ he thus got:

\begin{thm}[Zariski \cite{z}]
For any rational plane curve of even degree $2d-2$ with maximal number of cusps and nodes
$$
\pi_1(\PP^2- C)\quad\cong\quad\pi_1(\ufami_{1,d}).
$$
\end{thm}

Moreover he studies partial smoothings of such curves in which he allows cusps to resolve into
a node and a vertical tangent and nodes to resolve into two vertical tangents.
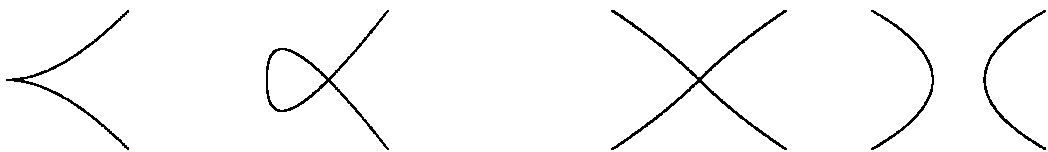
\unitlength=2.3mm
\begin{picture}(70,20)(-5,0)

\bezier{200}(0,10)(3,10)(7,14)
\bezier{200}(0,10)(3,10)(7,6)

\bezier{300}(15,10)(15,5)(22,14)
\bezier{300}(15,10)(15,15)(22,6)

\bezier{200}(35,14)(38,12)(40,10)
\bezier{200}(35,6)(38,8)(40,10)
\bezier{200}(40,10)(42,12)(45,14)
\bezier{200}(40,10)(42,8)(45,6)

\bezier{300}(50,14)(57,10)(50,6)
\bezier{300}(60,14)(53,10)(60,6)

\end{picture}

His result states, that in all but the nodal degree $6$ case, (d=4), such a deformation leads to a curve with complement
which has cyclic fundamental group $\ZZ_{2d-2}$.

In reverence to Zariski we
finish this paper with a proof of the following claim:

\begin{prop}
\labell{tribute}
For any plane curve, which admits a generic plane section of a discriminant $\dfami_{n,d}$ as above as limiting curve,
but has less singular points, the fundamental group of the complement is cyclic of degree
$(n+1)(d-1)^n$ except in the following cases:
\begin{enumerate}
\item
the exceptional case of Zariski, when a sextic has a generic plane section
of the discriminant $\dfami_{1,4}$ as a limiting curve with only
nodes smoothed. Then the fundamental group is isomorphic to $\psl$,
\item
when a duodectic has a generic plane section of the discriminant $\dfami_{2,3}$
as a limiting curve with only nodes smoothed. Then the
fundamental group of the complement is isomorphic to $\slz$.
\end{enumerate}
\end{prop}

\section{the discriminant polynomial}

The goal of this chapter is to gather some information about the polynomial $\dip$
which describes the discriminant $\dfami_{n,d}$ and its affine cone. 

Viewed as a polynomial in the distinguished parameter $u_{d,0,...,0}$ -- which now and hereafter
we will denote by $z$ -- we can immediately associate the \emph{bifurcation polynomial} $\bip=discr_z\dip$
and the leading coefficient $\lep$.
We collect some numerical invariants of the polynomials first to find then some
special features of their zero sets, which we are going to exploit in later sections.

\begin{lemma}
\labell{degp}
The discriminant polynomial $\dip$ is of degree
$$
\deg \dip\quad = \quad (n+1)(d-1)^n.
$$
\end{lemma}

\proof
This is well-know, but our proof is given here to get used to
the strategy employed in subsequent proofs.

Consider a generic pencil of degree $d$ hypersurfaces in $\pn$.
By definition of the discriminant the number of hypersurfaces
which are not smooth is $\deg\dip$. In fact they fail to be smooth in precisely one
ordinary double point by the genericity assumption.
All other hypersurfaces are smooth and the base locus of the pencil
is a smooth complete intersection of two smooth hypersurfaces.

To employ a topological Euler number calculation we recall that
\begin{enumerate}
\item
the Euler number of a smooth hypersurface of degree $d$ in $\pn$ is
$$
\euhyp\quad=\quad n+1+\frac{(1-d)^{n+1}-1}d
$$
\item
the Euler number of a hypersurface of degree $d$ in $\pn$ is
$\euhyp+(-1)^n$, if it is regular except for a single ordinary double point.
\item
the Euler number of a smooth complete intersection of two hypersurfaces
of degree $d$ in $\pn$ is 
$$
\euci=n+1+(1-d)^n(n-1)+2\frac{(1-d)^n-1}d.
$$
\end{enumerate}
Now we decompose $\pn$ into the base locus and its complement.
This complement is covered by the family of hypersurfaces each deprived
of the base locus. This family in turn can be decomposed into the family of smooth hypersurfaces
and the singular hypersurfaces.
So the Euler number of $\pn$ is obtained from the three corresponding summands:
\begin{eqnarray*}
& & e(\pn)\\
 & = & \euci + (e(\pe)-\deg\dip)(\euhyp-\euci)+\deg\dip(\euhyp-\euci+(-1)^n)\\
& = & 2\euhyp -\euci +(-1)^n\deg\dip
\end{eqnarray*}
Now all the numerical values are known except for $\deg\dip$, so we deduce
\begin{eqnarray*}
\deg\dip & = & (-1)^n\left(-e(\pn)+2\euhyp-\euci\right)\\
& = & (d-1)^n(n+1).\\[-15mm]
\end{eqnarray*}
\qed

\begin{lemma}
\labell{degzp}
The discriminant polynomial $\dip$ as a polynomial in the
coefficient $z$ of $x_0^d$ only has degree
$$
\deg_z\dip\quad=\quad(d-1)^n.
$$
\end{lemma}

\proof
The degree $\deg_z\dip$ is equal to the number of singular hypersurfaces
in a generic pencil of hypersurfaces of degree $d$ in $\pn$ containing
the multiple hyperplane $x_0^d$, but not counting this hyperplane.
Again we evaluate the topological Euler number of $\pn$ using the
decomposition into the hyperplane $x_0=0$ and its complement.

Note that the complement is covered by a family of hypersurfaces
parameterised by $\CC$, each deprived of the base locus which is
a smooth hypersurface of degree $d$ on the hyperplane $x_0=0$ which
is projective space of dimension $n-1$.
Moreover precisely $\deg_z\dip$ of these hypersurfaces are not smooth and
in fact regular except for a single ordinary double point.
The Euler number calculation yields the relation
\begin{eqnarray*}
&&e(\pn)\\
 & = & e(\pp^{n-1})+(e(\CC)-\deg_z\dip)\left(\euhyp-e_{n-1,d}\right)+\deg_z\dip(\euhyp-e_{n-1,d}+(-1)^n).
\end{eqnarray*}
Again we solve for the unknown with the numerical values provided above and get
\begin{eqnarray*}
\deg_z\dip & = &
(-1)^n\left(e(\pn) - e(\pp^{n-1}) - \euhyp + e_{n-1,d}\right)\\
& = & (d-1)^n\\[-15mm]
\end{eqnarray*}
\qed

\begin{lemma}
\labell{leadc}
The discriminant polynomial $\dip$ as a polynomial in the coefficient
$z$ of $x_0^d$ only has leading coefficient
$$
\lep\quad=\quad p_{n-1,d}^{d-1}.
$$
\end{lemma}

\proof
For any pencil of hypersurfaces of degree $d$ in $\pn$ containing the
multiple hyperplane $x_0^d$ we can compute the degree of the discriminant
polynomial in $z$ by an evaluation of topological Euler numbers again.
\begin{eqnarray*}
\deg_z & = & (-1)^n\left(e(\pn) - e(\pp^{n-1}) - e(\CC)(\euhyp - e(Bs))\right)\\
& = & (-1)^n\left(1- \euhyp + e(Bs)\right).
\end{eqnarray*}
where the base locus $Bs$ and $\deg_z$ depend on the pencil.

Hence the degree $\deg_z$ drops if and only if the Euler number of $Bs$
differs by a positive multiple of $(-1)^{n-1}$ as compared to $e_{n-1,d}$.

That change occurs if and only if $Bs$ is singular, which is a condition on the restriction
to the hyperplane $x_0=0$ of any general degree $d$ hypersurface of the pencil.
\\

We conclude that the degree $\deg_z$ drops if and only if the discriminant
polynomial in the appropriate variables -- the parameters of monomials
not containing $x_0$ -- vanishes.

The multiplicity of this factor is as claimed since degrees have to match:

The discriminant polynomial $\dip$ is homogeneous of degree $(n+1)(d-1)^n$ which hence
must be equal to the product of the degree $nk(d-1)^{n-1}+(d-1)^n$ of the leading coefficient
$p^k_{n-1,d}z^{(d-1)^n}$ and $\deg_z\dip$, so we infer our as $k=d-1$.
\qed

\begin{lemma}
\labell{copc}
The discriminant polynomial $\dip$ as a polynomial in the coefficient
$z$ of $x_0^d$ only has coprime coefficients.
\end{lemma}

\proof
By the preceding lemma
the leading coefficient has a unique irreducible factor $p_{n-1,d}$.
So the coefficients are not coprime only if $p_{n-1,d}$ is a 
factor of each.

In that case the zero set of $p_{n-1,d}$ belongs to the discriminant and must be
equal to the discriminant since both are irreducible.

This is not true because there are singular hypersurfaces which are regular
when restricted to $x_0=0$ so contrary to our assumption the coefficients are
coprime.
\qed

\begin{lemma}
\labell{degq}
The bifurcation polynomial $\bip$ is homogeneous of degree
$$
(2n+1)(d-1)^n\left((d-1)^n-1\right).
$$
\end{lemma}

\proof
The polynomial $\bip$ is obtained as the discriminant of $\dip$ with
respect to the variable $z$ and so is homogeneous itself. In fact it can be computed from the
Sylvester matrix of $\dip$ and $\del_z\dip$; up to a factor consisting
of the leading coefficient polynomial $\lep$ it is the determinant of that
matrix.

Hence it is sufficient to add the degrees along the diagonal for any reordering of the matrix.
In fact we can arrange on this diagonal $\deg_z\dip$ times the leading coefficient of degree
$\deg\dip-\degz\dip$ and $\degz\dip-1$ times the constant coefficient of degree $\deg\dip$.

By the above we have to subtract the degree of the leading coefficient and thus we get
\begin{eqnarray*}
\deg\bip & = & \degz\dip(\deg\dip-\degz\dip)+(\degz\dip-1)\deg\dip\\
&&
-(\deg\dip-\degz\dip)\\
&=&
(\degz\dip-1)(\deg\dip-\degz\dip)+(\degz\dip-1)\deg\dip\\
&=&
(\degz\dip-1)(2\deg\dip-\degz\dip)\\
&=& \left((d-1)^n-1\right)\left(2(n+1)(d-1)^n-(d-1)^n\right)\\
&=&
(2n+1)\left((d-1)^n-1\right)(d-1)^n
\end{eqnarray*}
\qed

In the next instances we consider $\dip$ to be weighted homogeneous. The weight we assign to
each $u_\nu$ is equal to the exponent of $x_0$ in the monomial of which it is the parameter:
$$
wt(u_\nu)\quad = \quad \nu_0.
$$
In particular the weight is zero if $x_0$ does not occur and it is $d$ precisely 
in the case of the parameter $z$.
\\

Another shorthand will then be $v_\ncount$ for the parametrical coefficients $u_{\nu_{(\ncount)}}$
of the mono\-mial $x_0^{d-1}x_\ncount$, so $\nu_{(\ncount)}=(\nu_0,...,\nu_n)$ with no non-zero
components except for $\nu_0=d-1$ and $\nu_\ncount=1$.
We call the $v_\ncount$ the linear parametrical coefficients.

\begin{lemma}
\labell{wdegp}
$\dip$ is weighted homogeneous of degree
$$
\wdeg\dip\quad=\quad d(d-1)^n.
$$
\end{lemma}

\proof
The leading term of $\dip$ is of degree $(d-1)^n$ in $z$ which is of weight $d$
with coefficient in the coefficients of the monomials not involving $x_0$ at all.
So all variables in the leading coefficient are of weight zero and the claim follows.
\qed

\begin{lemma}
\labell{wdegq}
$\bip$ is weighted homogeneous of degree
$$
\wdeg\bip\quad=\quad d(d-1)^n\left((d-1)^n-1\right).
$$
\end{lemma}

\proof
The discriminant of $\dip$ with respect to $z$ is $\bip$, so
\begin{eqnarray*}
\wdeg\bip & = &
\wdeg z\degz\dip(\degz\dip-1)
\end{eqnarray*}
and the numerical values given in the lemmas above yield the claim.
\qed

\begin{lemma}
\labell{degvq}
$\bip$ as a polynomial in the linear coefficients $v_\ncount$ is of degree
$$
\degv\bip\quad=\quad d(d-1)^{n-1}\left((d-1)^n-1\right).
$$
\end{lemma}

\proof
The linear coefficients are of weight $d-1$, so we get the upper bound for
$\degv\bip$ to be $\wdeg\bip/(d-1)$.

The existence of at least one non-trivial coefficient
can be deduced from the special family
$$
\sum a_\ncount x_\ncount ^d+\sum v_\ncount x_\ncount x_0^{d-1}+zx_0^d.
$$
For all $a_\ncount =1$ and all $v_\ncount $ positive real of sufficiently distinct magnitude
$$0<v_1<<v_2<<\cdots<<v_n$$ the discriminant polynomial has simple roots only
so the bifurcation polynomial for the family is non-zero.

On the other hand the bifurcation polynomial of our special family is 
weighted homogeneous again with $\wdeg=d(d-1)^n((d-1)^n-1)$ to which in fact
only the $v_\ncount $ contribute since all $a_\ncount$ have weight $0$. So from the non-triviality
above we conclude our claim.
\qed

\begin{lemma}
\labell{degcv}
Consider $\bip$ as a polynomial in the parameters $v_\ncount $ with coefficients in $\CC[u_\nu|\nu_0<d-1]$.
Given a monomial of degree $\degv\bip$ in the parameters $v_\ncount $, its coefficient $c$
in $\bip$ is a polynomial in the parameters $u_\nu'$ of monomials $x^\nu$
not containing $x_0$ and it has degree
$$
\deg c \quad = \quad \left(2n(d-1)-1\right)(d-1)^{n-1}\left((d-1)^n-1\right).
$$
(Of course the coefficient may be zero.)
\end{lemma}

\proof
The degree is just the difference between $\deg\bip$ and $\degv\bip$.
The other claim follows from the fact that the leading coefficient has weighted degree $0$, so must
be a polynomial in the weight $0$ parameters $u_\nu'$.
\qed

\begin{lemma}
\labell{coprimv}
Consider $\bip$ as a polynomial in $(\CC[u_\nu|\nu_0<d-1])[v_\ncount ]$.
Then there is no nontrivial common divisor of all its coefficient polynomials in
$\CC[u_\nu|\nu_0<d-1]$.
\end{lemma}

\proof
By construction a polynomial $f$ belongs to the zero set of $\bip$ if and only if
at least one of the polynomials $f+zx_0^d, z\in\CC$ has more then a single ordinary
double point singularity
or the restriction $f|_{x_0=0}$ has more than
a single ordinary double point singularity.

If $f$ is any polynomial then some perturbation $\tilde f$
of $f$ by terms $v_\ncount x_\ncount x_0^{d-1}$ has the property
that $\tilde f$ has non-degenerate critical points only.
Moreover by changing the perturbation ever so slightly
we may assume, that $\tilde f$ has even no multiple critical
values. Hence $\tilde f$ belongs to the zero set of $\bip$ if
and only if the restriction $\tilde f|_{x_0=0}$ is non-generically singular.

The zero set of a common factor of all coefficients is
either a divisor or empty since $\bip$ is non-trivial.

Moreover by the preceding lemma the coefficients
of monomials in $\CC[v_\ncount ]$ of highest degree contain only
parameters $u'$, hence the same must be true for a common
factor of all coefficients.

Therefore if a polynomial $f$ belongs to the zero set of a common factor
then so does every perturbation $\tilde f$ as above.

Hence a polynomial can only belong to the zero set
if its restriction $f|_{x_0=0}$ is non-generically singular.
Since the set of polynomials with non-generically singular restriction
to $x_0=0$ is of codimension two, the zero set of any
common factor is empty and therefore any such common
factor is a non-zero constant.
\qed

\section{Zariski arguments}

The ideas of Zariski as elaborated by Bessis \cite{bessis} in the affine set up
provide the tool to get hold of a presentation for the fundamental group
of the complement of the discriminant cone.

Most of the geometry involved relies on the notion of \emph{geometric element},
which refers to elements in a fundamental group which are represented by
paths isotopic to a boundary of a small disc transversal to a divisor.

In case of a punctured disc or affine line a basis for the fundamental group is called a \emph{geometric basis}
if it consists of geometric elements simultaneously represented by paths only
meeting in the base point.

Basic to our argument is the following observation made explicit by Bessis \cite{bessis}:

\begin{lemma}[Bessis]
\labell{bess}
For any affine divisors $D,E$ without common component there is an exact sequence of fundamental groups:
$$
\pi_1(\CC^N-D-E)\surj\pi_1(\CC^N-D)\tto1.
$$
with kernel normally generated by the geometric elements associated to $E$.
\end{lemma}

We denote now by $\CC^N$ the affine parameter space containing divisors
$\dfami,\hat\afami,\hat\bfami$, where the discriminant cone $\dfami$
is given by $\dip$, $\hat\afami,\hat\bfami$ are given by $p_{n-1,d}$ and $\bip$
respectively. 
Note that for notational convenience we supress the dependence on integers $n,d$ occasionally.

We project $\CC^N$ to $\CC^{N-1}$ along the distinguished parameter $z$ and
get divisors $\afami,\bfami$ defined by $p_{n-1,d}$ and $\bip$ again.
Of course by construction $\hat\afami,\hat\bfami$ are the pull backs of $\afami,\bfami$
along the projection and $\dfami$ is finite over $\CC^{N-1}$, branched along $\bfami$
and has $\hat\afami$ as its vertical asymptotes, hence $\afami$ is the pole locus
of $\dfami$ with respect to the projection.

\begin{lemma}
\labell{sabd} 
Suppose $L$ is a fibre of the projection such that its intersection $\dfami_L$
with the discriminant $\dfami$ consists of $\deg_z\dip$ points,
then there is a split exact sequence
$$
1\to\pi_1(L-\dfami_L)\to\pi_1(\CC^N-\hat\afami-\hat\bfami-\dfami)
\to
\pi_1(\CC^{N-1}-\afami-\bfami)\to1.
$$
with a splitting map which takes geometric elements associated to $\bfami$
to geometric elements associated to $\hat\bfami$.
\end{lemma}

\proof
In fact over the complement of $\afami\cup\bfami$ the discriminant is a finite topological cover
and its complement is a locally trivial fibre bundle with fibre the affine line punctured at
$\deg_z\dip$ points.
The exact sequence is now obtained from the long exact sequence of that fibre bundle.
Exactness on the left follows from the fact that no free group of rank more than $1$ admits
a normal abelian subgroup.

The splitting map is induced by a topological section given by 
a real number which bounds the modulus of all zeroes of the
given polynomial. Of course this bound can be chosen continuously in terms
of the coefficients of the polynomial.

The final observation is, that this topological section maps boundaries of small discs transversal
to $\bfami$ to boundaries of small discs transversal to $\hat\bfami$ and disjoint to any other
divisor.
\qed

We should nevertheless note, that boundaries of arbitrarily small discs transversal to $\afami$ are mapped
to boundaries of discs transversal to $\hat\afami$ but also intersecting $\dfami$.
\\

The group in the middle is hence determined as the semi-direct product of the other two
by a map of $\pi_1(\CC^{N-1}-\afami-\bfami)$ to the automorphism group of $\pi_1(L-\dfami_L)$.

This has an immediate corollary on the level of presentations:

\begin{lemma}
\labell{pabd}  
Suppose there is a presentation for the fundamental group of the base
$$
\pi_1(\CC^{N-1}-\afami-\bfami)\quad\cong\quad
\langle r_\a | \rfami_q\rangle.
$$
in terms of geometric generators then there is a presentation
$$
\pi_1(\CC^n-\hat\afami-\hat\bfami-\dfami)
\quad\cong\quad
\langle t_i,\hat r_\a | \hat r_\a t_i\inv\hat r_\a\inv\phi_\a(t_i),\rfami_q\rangle.
$$
where $\phi_\a$ is the automorphism associated to $r_\a$ and $t_i$ is a free
geometric basis for a generic vertical line $L$ punctured at $L\cap\dfami$.
\end{lemma}

\newcommand{\vplus}{v_\Sigma} 

To investigate possible presentation for the fundamental group of the base, we want to
exploit another projection $\CC^{N-1}$ to $\CC^{N-2}$ along some linear combination $\vplus$
of the parameters $v_\ncount$.

We note, that $\afami$ is then the pull back of its image $\bar\afami$
since the equation of $\afami$ does not contain any of the $v_\ncount$.
On the other hand $\bfami$ is more sensitive to the choice of $\vplus$.
But without getting explicit we can prove existence of a suitable projection:

\begin{lemma}
\labell{projv}
\begin{sloppypar}
There is a linear combination $\vplus$ of the $v_\ncount$ such that the projection \mbox{$p_v:\CC^{N-1}\to\CC^{N-2}$}
along $\vplus$ has the following property:
\begin{quote}
There exists a divisor $\bar\cfami$ such that no component of its pull-back $\cfami$ to $\CC^{N-1}$
is a component of $\bfami$ and such that the induced map
$p_v|:\bfami\to \CC^{N-2}$ is a topological finite cover over the complement of $\bar\cfami$.
\end{quote}
\end{sloppypar}
\end{lemma}

\proof
The set of singular values for the induced map $p_v|:\bfami\to \CC^{N-2}$ is a divisor $\bar\cfami$, since
we equip $\bfami$ with its reduced structure. Of course in its complement $\bfami$ must be a topological
fibration hence a topological covering.

A common component of $\bfami$ and the divisor $\cfami$ can be detected
algebraically as a nontrivial factor
of $\bip$ which is independent of the variable $\vplus$.
It thus suffices to show that for a suitable choice of $\vplus$ there is no such factor.
\\

We decompose the polynomial algebra $\CC[u_\nu]$ according to the degree $\deg_v$
of each monomial considered as a monomial in the $v_\ncount$ only.
With respect to the $\deg_v$ decomposition we have the summand $q_{max}$ of $\bip$ of
highest degree. It is of degree $d(d-1)^{n-1}((d-1)^n-1)$ by lemma \reff{degvq} and
its coefficients are in $\CC[u'_\nubold]$ by lemma \reff{degcv}.

Therefore $q_{max}$ defines a proper hypersurface in some trivial affine bundle over $\PP^{n-1}$.
If we replace the projective coordinates $v_1:v_2:...:v_n$ of $\PP^{n-1}$ by new ones
$v_\ncount'$ in such a way that a point in the complement
has coordinates $(0:...:0:1)$ then $\bip$ is a polynomial of degree $d(d-1)^{n-1}((d-1)^n-1)$
in the variable $\vplus=v_n'$ with leading coefficient in $\CC[u_\nu']$.

With that choice a non-trivial factor of $\bip$ independent of $\vplus$ may only depend on the $u_\nu'$.
\\

But we know already from the proof of lemma \reff{coprimv}, that no
divisor defined in terms of the $u'_\nubold$ only can be a component of
$\bfami$, hence there is no common component of $\cfami$ and $\bfami$.
\qed

We suppose from now on a projection $p_v:\CC^{N-1}\to\CC^{N-2}$ as in the
lemma has been fixed with $L'$ a generic fibre.

\begin{lemma}
\labell{genab}  
Suppose there are geometric elements $r_a$ associated to $\afami$
and a geometric basis consisting of elements $r_b$ of $\pi_1(L'-\bfami_{L'})$
such that the $r_a$ generate $\pi_1(\CC^{N-1}-\afami)$, 
then the $r_a$ and $r_b$ together generate
$$
\pi_1(\CC^{N-1}-\afami-\bfami).
$$
\end{lemma}

\proof
From the lemma of Bessis we infer that given the elements $r_a$
which generate $\pi_1(\CC^{N-1}-\afami)$ there are geometric elements $r_c$
associated to $\cfami$ which can be
taken in the complement of $\bfami$ such that together they generate
$$
\pi_1(\CC^{N-1}-\afami-\cfami)\quad\cong\quad\pi_1(\CC^{N-2}-\bar\afami-\bar\cfami).
$$
We are in fact in a situation similar to that of lemma \reff{sabd} since $\CC^{N-1}-\afami-\bfami-\cfami$
fibres locally trivial over $\CC^{N-2}-\bar\afami-\bar\cfami$ with fibre
the appropriately punctured complex line $L'-\bfami_{L'}$.
In the associated exact sequence
$$
\pi_1(L_b-\bfami_L)\tto\pi_1(\CC^{N-1}-\afami-\bfami-\cfami)\tto\pi_1(\CC^{N-2}-\bar\afami-\bar\cfami)
$$
the elements $r_b$ generate the group on the left hand side and the images
of the $r_a$ and $r_c$ generate the group on the right hand side, so together they
generate the group in the middle.

But then we can apply the result of Bessis once more to see, that together they generate
$$
\pi_1(\CC^{N-1}-\afami-\bfami)
$$
and that the elements $r_c$ are trivial in that group, so our claim holds.
\qed

\begin{lemma}
\labell{pad}  
Suppose there is a presentation for the fundamental group of $\CC^{N-1}-\afami$
$$
\pi_1(\CC^{N-1}-\afami)\quad\cong\quad
\langle r_a | \rfami_a\rangle.
$$
in terms of geometric generators and that $\pi_1(L_b-\bfami_L)$ is generated
by a geometric basis $r_b$ then there is a presentation
$$
\pi_1(\CC^n-\hat\afami-\dfami)
\quad\cong\quad
\langle t_i,\hat r_a | t_i\inv\phi_b(t_i),\hat r_a t_i\inv\hat r_a\inv\phi_a(t_i),\rfami_a\rangle.
$$
where $\phi_a$ ($\phi_b$) is the automorphism associated to $r_a$ ($r_b$), $t_i$ is a free
geometric basis of $\pi_1(L-\dfami_L)$ and the $\hat r_a$ are lifts of $r_a$ by the topological section.
\end{lemma}

\proof
We apply the lemma of Bessis to the result of Lemma \reff{pabd} . In the presentation we have simply
to set the geometric generators associated to $\hat\bfami$ to be trivial and to discard them from
the set of generators.
\qed

\begin{lemma}
\labell{pd}  
The fundamental group of $\CC^N-\dfami$ has a presentation
$$
\pi_1\quad\cong\quad
\langle t_i | t_i\inv\phi_b(t_i),\rho_at_i\inv \rho_a\inv\phi_a(t_i)\rangle.
$$
where
\begin{enumerate}
\item
the $t_i$ form a geometric basis of $\pi_1(L-\dfami_L)$,
\item
the $\rho_a$ can be expressed in terms of $t_i$ such that $\hat r_a\rho_a\inv$ is a geometric element
associated to $\hat\afami$ and a lift of $r_a$,
\item
the $\phi_a$ are the braid monodromies associated to $r_a$.
\end{enumerate}
\end{lemma}

\proof
Since each $\hat r_a$ arises from an element transversal to $\hat\afami$ it is equal to a geometric element $\hat r_a'$ for
$\hat\afami$ up to some $\hat\rho_a$ expressible in terms of geometric elements for $\dfami$.

But in the complement of $\dfami$ only, the elements $\hat r_a$ become trivial and their set relations $\rfami_a$
may be discarded, so from the previous lemma we get to our claim.
\qed

The claim of the lemma is of course only an intermediate step on our way to give a presentation
of the fundamental group.
Obviously we have to make the relations explicit in the sense that every relation is given in terms
of the chosen generators only.

\paragraph{Remark}
We are very lax about the base points. They should be chosen in such a way
that all maps of topological spaces are in fact maps of pointed spaces.
(In particular in the presence of a topological section there is no choice left;
in the fibre and in the total space the base point is the intersection of the
section with the fibre and its projection to the base yields the base point there.)

\section{Brieskorn Pham unfolding}
\label{affine}

In this section our aim is twofold, first to propose a distinguished set of generators for $\pi_1(\CC^N-\dfami)$
and second to give explicitly an exhaustive set of relations associated to a geometric basis for the complement
of $\bfami$.

So first we pick some distinguished fibres $L_v$ of the projection
$p_z:\CC^N\to\CC^{N-1}$ along the variable $z$ where in each case $L_v-\dfami_L$ can 
be equipped with a distinguished geometric basis by the method of Hefez and Lazzeri \cite{hl}.
For later use in section \reff{projective} we establish a relation between different such bases.

For our second aim we exploit the relation with the versal unfolding of
Brieskorn Pham singularities to make those relations
explicit which arise from geometric generators associated to $\bfami$.

\subsection{Hefez Lazzeri path system}
\label{hl-path}

First we want to describe a natural geometric basis for some fibres of the projection
$p:\CC^N\to\CC^{N-1}$.
Since we follow Hefez and Lazzeri \cite{hl} we will call such bases accordingly.
We note first that fibres $L_u$ of the projection correspond to affine pencil of polynomials
$$
f_u(x_1,...,x_n)-z
$$
and their discriminant points $\dfami_L$ are exactly the $z$ such that the $z$-level
of $f$ is singular.
As in \cite{hl} we restrict our attention to the linearly perturbed Fermat polynomial:
$$
f\quad=\quad\sum_{\ncount=1}^n (x_\ncount^d-dv_\ncount x_\ncount).
$$

In that family the discriminant points for any generic pencil are in bijection
to the elements in the multiindex set of cardinality $(d-1)^n$:
$$
\ind\quad=\quad\{\,(i_1,...,i_n)\,|\,1\leq i_\nu\leq d-1\,\}.
$$
More precisely we get an expression for the critical values from \cite{hl}:

\begin{lemma}[Hefez Lazzeri]
\labell{laz-dis}
The polynomial defining the critical value divisor is given by the expansion
of the formal product
$$
\prod_{\ibold \in\ind}
\left(-z+(d-1)\sum_{\ncount=1}^n\xi^{i_\ncount}v_\ncount^{\frac{d}{d-1}}\right).
$$
\end{lemma}

As an immediate corollary we show that discriminant sets are equal for suitably
related parameter values:

\begin{lemma}
\labell{twist}
The discriminant of the linearly perturbed Fermat polynomial
is invariant under the multiplication of any $v_\ncount$ by a
$d$-th root of unity.
\end{lemma}

\proof
From the expansion above we see that the discriminant polynomial is a polynomial in
$v_\ncount^{\frac{d}{d-1}}$ but of course it is also a polynomial in $v_\ncount$, hence it must
be a polynomial in $v^d_\ncount$ since that is the least common power of both.
Then it is obviously invariant under multiplying $v_\ncount$ by a $d$-th root.
\qed

We also note the following induction property of the discriminant points:

\begin{lemma}
\labell{circ}
The critical values of $f$ are distributed on circles of radius $(d-1)|v_\ncount|^{\frac{d}{d-1}}$
centred around the critical values of the polynomial
$$
f'\quad=\quad \sum_{\ncount=1}^{n-1}(x_\ncount^d-dv_\ncount x_\ncount).
$$
\end{lemma}

\proof
Again we can use lemma \reff{laz-dis}.
A formal zero of the discriminant polynomial for $f$ differs by a term
$(d-1)v_\ncount^{\frac{d}{d-1}}$ from a zero of the discriminant polynomial of $f'$, and that
difference is of the claimed modulus.
\qed

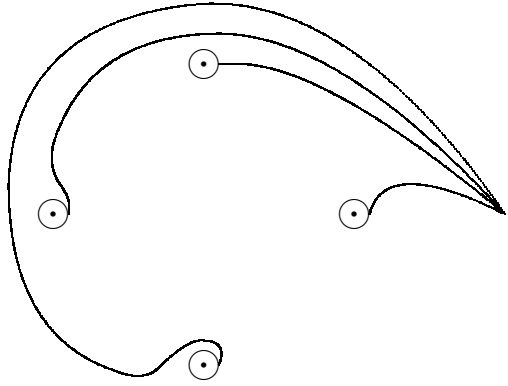
\begin{figure}[ht]
\label{base}
\begin{center}
\setlength{\unitlength}{.4mm}
\begin{picture}(180,150)(-70,-70)

\put(50,0){\circle{10}}
\put(0,50){\circle{10}}
\put(0,-50){\circle{10}}
\put(-50,0){\circle{10}}
\put(50,0){\circle*{2}}
\put(0,50){\circle*{2}}
\put(0,-50){\circle*{2}}
\put(-50,0){\circle*{2}}


\bezier{130}(55,0)(60,20)(100,0)


\bezier{30}(5,50)(8,50)(11,50)
\bezier{300}(11,50)(40,50)(100,0)

\bezier{300}(-50,24)(-52,16)(-48,10)
\bezier{300}(-48,10)(-44,5)(-45,0)
\bezier{300}(-50,24)(-38,60)(5,60)
\bezier{300}(5,60)(45,60)(100,0)

\bezier{80}(5,-50)(8,-42)(-1,-42)
\bezier{80}(-1,-42)(-6,-42)(-15,-51)

\bezier{80}(-15,-51)(-20,-56)(-30,-52)
\bezier{200}(-30,-52)(-65,-40)(-65,10)
\bezier{200}(-65,10)(-63,65)(00,70)
\bezier{200}(0,70)(50,72)(100,0)
\end{picture}
\caption{Hefez Lazzeri system in case $n=1$, $d-1=4$}
\end{center}
\end{figure}

We assume now that all $v_\ncount$ are positive real and of sufficiently distinct modulus
$$
|v_n|\ll\cdots\ll|v_1|.
$$
In case $n=1$ we define the Hefez Lazzeri geometric basis as indicated in figure \reff{base} for $d-1=4$,
where each geometric generator is depicted as a tail and a loop around a critical value.

Of course the geometric element associated to a loop-tail pair is represented by a closed
path based at the free end of the tail which proceeds along the tail, counterclockwise around
the loop and back along the tail again.

The $d-1$ elements of the base are denoted by $t_1,...,t_{d-1}$, such that the corresponding critical
values are enumerated counterclockwise starting on the positive real line.

Each tail is thus also labeled by a number in $\{1,...,d-1\}$.

For the inductive step we suppose that the elements of a Hefez Lazzeri base for
$$
f'\quad=\quad \sum_{\ncount=1}^{n-1}(x_\ncount^d-dv_\ncount x_\ncount).
$$
are given by tail loop pairs each labeled by some multi-index in $I_{n-1,d}$.
By assumption $v_n$ is sufficiently small compared to $v_{n-1}$ so we may assume that
all critical values of
$$
f \quad=\quad\sum_\ncount (x_\ncount^d-dv_\ncount x_\ncount).
$$
are inside the loops and in fact distributed at distance $(d-1)|v_\ncount|^{\frac{d}{d-1}}$
from their centres.

In the inductive step each loop and its interior are erased and replaced by scaled copy
of the Hefez Lazzeri base for the $n=1$ case.
Each tail-loop pair with label $i_n\in\{1,...,d-1\}$ in an inserted disc fits with a tail labeled
by some $\ibold'=i_1i_2\cdots i_{n-1}$ to form a tail-loop pair representing an element of the 
new Hefez-Lazerri base which is labeled by $\ibold =i_1i_2\cdots i_{n-1}i_n$.

\begin{figure}[ht]
\setlength{\unitlength}{.4mm}
\begin{picture}(180,150)(-250,-70)

\put(50,0){\circle*{10}}
\put(0,50){\circle*{10}}
\put(0,-50){\circle*{10}}
\put(-50,0){\circle*{10}}


\bezier{130}(55,0)(60,20)(100,0)


\bezier{30}(5,50)(8,50)(11,50)
\bezier{300}(11,50)(40,50)(100,0)

\bezier{300}(-50,24)(-52,16)(-48,10)
\bezier{300}(-48,10)(-44,5)(-45,0)
\bezier{300}(-50,24)(-38,60)(5,60)
\bezier{300}(5,60)(45,60)(100,0)

\bezier{80}(5,-50)(8,-42)(-1,-42)
\bezier{80}(-1,-42)(-6,-42)(-15,-51)

\bezier{80}(-15,-51)(-20,-56)(-30,-52)
\bezier{200}(-30,-52)(-65,-40)(-65,10)
\bezier{200}(-65,10)(-63,65)(00,70)
\bezier{200}(0,70)(50,72)(100,0)

\put(-170,0){
\unitlength=.28mm
\begin{picture}(180,150)

\put(50,0){\circle{10}}
\put(0,50){\circle{10}}
\put(0,-50){\circle{10}}
\put(-50,0){\circle{10}}
\put(50,0){\circle*{2}}
\put(0,50){\circle*{2}}
\put(0,-50){\circle*{2}}
\put(-50,0){\circle*{2}}

\bezier{130}(55,0)(60,20)(100,0)

\bezier{30}(5,50)(8,50)(11,50)
\bezier{300}(11,50)(40,50)(100,0)

\bezier{300}(-50,24)(-52,16)(-48,10)
\bezier{300}(-48,10)(-44,5)(-45,0)
\bezier{300}(-50,24)(-38,60)(5,60)
\bezier{300}(5,60)(45,60)(100,0)

\bezier{80}(5,-50)(8,-42)(-1,-42)
\bezier{80}(-1,-42)(-6,-42)(-15,-51)

\bezier{80}(-15,-51)(-20,-56)(-30,-52)
\bezier{200}(-30,-52)(-65,-40)(-65,10)
\bezier{200}(-65,10)(-63,65)(00,70)
\bezier{200}(0,70)(50,72)(100,0)
\end{picture}
\put(-175,0){\setlength{\unitlength}{5.4mm}\begin{picture}(10,10)(5,5)
\bezier{30}(0,5)(0,7.5)(2,9)
\bezier{30}(2,9)(5,11)(8,9)
\bezier{30}(8,9)(10,7.5)(10,5)
\bezier{30}(10,5)(10,2.5)(8,1)
\bezier{30}(8,1)(5,-1)(2,1)
\bezier{30}(2,1)(0,2.5)(0,5)

\bezier{30}(8,9)(11,7)(13.4,5.4)
\bezier{30}(8,1)(11,3)(13.4,4.6)

\end{picture}}
}
\end{picture}
\caption{Hefez Lazzeri system in case $n=2$, $d-1=4$}
\end{figure}
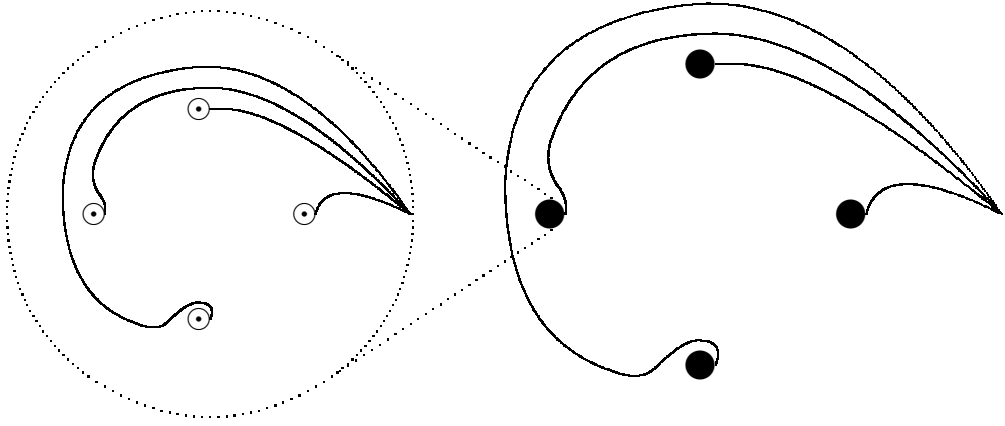

The element $\d_0$ represented by a path enclosing all critical values counter-clockwise
deserves our special attention.
In the case $n=1$ we see immediately that as an element in the fundamental group it can
be expressed as
$$
t_{d-1}t_{d-2}\cdots t_2t_1.
$$

To exploit the inductive construction for the general case we define $t_{\ibold'}^+$,
$\ibold'\in I_{n-1,d}$, in the Hefez Lazzeri fibre of 
$$
f \quad=\quad\sum_\ncount (x_\ncount^d-dv_\ncount x_\ncount)
$$
to be represented by the same tail-loop pair as the element $t_{\ibold'}$ in
the Hefez Lazzeri fibre of
$$
f'\quad=\quad \sum_{\ncount=1}^{n-1}(x_\ncount^d-dv_\ncount x_\ncount).
$$

\begin{lemma}
\labell{bundle}
For any $\ibold'\in I_{n-1,d}$ there is a relation
$$
t_{\ibold'}^+\quad=\quad t_{\ibold'(d-1)}t_{\ibold'(d-2)}\cdots t_{\ibold'2}t_{\ibold'1},
$$
where ${\ibold'(d-1)},{\ibold'(d-2)},\ldots, {\ibold'2},{\ibold'1}$ are the obvious elements of $\ind$.
\end{lemma}

\proof
The loop of $t_{\ibold'}^+$ is the path which encloses counter-clockwise $d-1$ of the
critical values which are inserted in the inductive step.
So it is homotopic to the product with descending indices of the loop-tail pairs inserted into that loop.

This relation is preserved under appending the tail of $t_{\ibold'}$ and so the claim follows.
\qed

To formulate the general claim we define a total order $\prec_0$ on $\ind$ to be
the lexicographical order with respect to the total order $>$ (!) on each
component.
Accordingly we introduce an order preserving enumeration function
$$
\ibold_0:\,(\{1,...,(d-1)^n\},<)\,\tto\,(\ind,\prec_0).
$$

\begin{lemma}
\labell{cox}
Suppose the element $\d_0$ is represented by a path which
encloses counter-clockwise
all critical values of a Hefez Lazzeri fibre, then
$$
\d_0\quad=\quad\prod_{k=1}^{(d-1)^n}t_{\ibold_0(k)}.
$$
\end{lemma}

\proof
In case $n=1$ we have given above an expression for $\d_0$ which is the expression claimed
here, as we have taken care that the order given by $\prec_0$ on $I_{1,d}$ is
$$
d-1 \,\prec_0\, d-2 \,\prec_0\, \cdots\,\prec_0\, 2 \,\prec_0\, 1.
$$
Suppose now the claim is true for $n-1$, that is
$$
\d_0'\quad=\quad
\prod_{k=1}^{(d-1)^{n-1}}t_{\ibold_o'(k)}.
$$
From the definition of $t_{\ibold'}^+$ we deduce immediately
$$ 
\d_0\quad=\quad
\prod_{k=1}^{(d-1)^{n-1}}t_{\ibold_o'(k)}^+
$$
where we can replace each $t_{\ibold'}^+$ using the preceding lemma \reff{bundle}
which yields our claims as we have defined $\ibold_0$ in the appropriate way.
\qed

\subsection{paths and homotopies}
\labell{pathsandhomotopies}

As we noticed in lemma \reff{twist} the set of singular values
remains unchanged upon multiplication of any of the real
$v_\ncount$ by a $d$-th root of unity.
Hence we may employ the same Hefez-Lazzeri system of paths.

For any $\jbold\in\ind$ and the primitive $d$-th root of unity $\xi$ of least argument
we introduce the notation $t_\ibold(\jbold)$ for the elements of the Hefez-Lazzeri basis
in the fibre over the parameter $v(\jbold)=(v_1\xi^{1-j_1},...,v_n\xi^{1-j_n})$.

Of course the different fundamental groups $\pi_1(L_{v(\jbold)}-\dfami_{L_\jbold},(v(\jbold),z_0))
$ are related by isomorphisms provided by
paths between their base points. Let $\oo_\jbold$ denote the path
$$
s\,\mapsto\,(v_1\xi^{s(1-j_1)},...,v_n\xi^{s(1-j_n)},z_0),
$$
then the corresponding isomorphism can be used to compare geometric generators in different fibres:

\begin{lemma}
\labell{hotopy}
Conjugation by a path $\oo_\jbold$ induces an isomorphism
$$
\oo^*_\jbold:\pi_1(L_{v(\jbold)}-\dfami_{L_\jbold},(v(\jbold),z_0))\to\pi_1(L_v-\dfami_L,(v,z_0))
$$
such that for $\ibold_0(1)=11\cdots 1\in\ind$ 
$$
t_\jbold
\quad=\quad
\oo^*_\jbold(t_{\ibold_0(1)}(\jbold))
$$
\end{lemma}

\proof
We consider case $n=1$ first. Consecutive critical values form an angle $\varth_d=\frac{2\pi}{d-1}$.
Let $\Phi_j(s)$ be an isotopy of the plane which fixes all points of modulus at least $z_0$ and rotates
rigidly all point of modulus at most that of the critical values till it is rotated by $(1-j)\varth$.

Then the homotopy class of $t_1$ in $L_v$ is mapped to the homotopy class of $t_j$ in $L_{v(j)}=L_v$.
Moreover we get from lemma \ref{laz-dis} that $\Phi_j$ preserves the discriminant along $\oo_j$,
i.e.\ $\Phi_j(s)$ maps the critical points in $L_v$ to the critical points in the fibre containing $\oo_j(s)$.
\\

If $t_1$ denotes the map which yields the closed path $t_1$ we can define a map
$$
s\,\mapsto\, (p(\oo_j(s)),\Phi_j(s)[t_1])
$$
as a one parameter family of paths.

By the last observation above its image is disjoint from the discriminant, hence it is a homotopy in the
complement.

Note that for $s=0$ we get the path $t_1$. Proceeding along the other part of the boundary we follow
first $\oo_j$ then $t_j(j)$ and the $\oo_j$ backwards, hence the claim follows for $n=1$.
\\

In the general case the isotopy above has to be replaced by a 'carousel' isotopy $\Phi_\jbold^{(n)}$ which can
be constructed inductively.

So we suppose that we have an isotopy $\Phi_{\jbold'}^{(n-1)}$, such that with
respect to the critical points of $f'$ and the path $\oo_{\jbold'}$
\begin{enumerate}
\item
$\Phi_{\jbold'}^{(n-1)}$ preserves the discriminant along $\oo_{\jbold'}$,
\item
$\Phi_{\jbold'}^{(n-1)}$ maps the homotopy class of $t_{\ibold'_o(1)}$
to that of $t_{\jbold'}$
\end{enumerate}

The critical points of $f$ are equi-distributed on small circles of fixed radius, cf.\ lemma \ref{circ},
around those of $f'$ along all paths $oo_\jbold$.

Now $\Phi_{\jbold'}^{(n-1)}$ can be slightly deformed without loosing its properties above
to get an isotopy $\Phi_{\jbold'1}^{(n)}$ which transforms each disc bounded by the above circles
by an affine translation.

It is readily checked that $\Phi_{\jbold'1}^{(n)}$ has the necessary properties with
respect to the critical points of $f$ and the path $\oo_{\jbold'1}$.
More generally  we may define $\Phi_{\jbold}^{(n)}$ to be $\Phi_{\jbold'1}^{(n)}$ composed
with $\Phi_{j_n}^{(1)}$ on each of the discs above.
\\

The claim now follows from the two properties as in the $n=1$ case;
the isotopy $\Phi_{\jbold}^{(n)}$ realises the claimed homotopy equivalence.
\qed

\subsection{Brieskorn-Pham monodromy}

The next aim is to determine the set of relations imposed on Hefez-Lazzeri generators
by the geometric generators associated to $\bfami$.

Our strategy is to show the genericity relative to the
bifurcation set $\bfami$ of a line $L'$ in the unfolding of the
Brieskorn Pham singularity
$$
x_1^d+x_2^d+\cdots+ x_n^d,
$$
by non-constant monomials of degree less than $d$, the \emph{truncated subdiagonal unfolding}.
\\

This can only be true if the corresponding pencil of polynomials degenerates only in
the most modest ways, i.e.\ all polynomials are Morse except for a finite set of polynomials which have
precisely one critical value less, due either to the presence of a degenerate critical point, necessarily of type $A_2$,
or to the coincidence of the values of two distinct non-degenerate critical points.

This necessary condition is actually met in our situation:

\begin{lemma}
\labell{a2gen}
There is a pencil of subdiagonally perturbed polynomials, each of which has only non-degenerate critical point,
except for a finite number of polynomials with a solitary degenerate critical point of
type $A_2$ of singular value distinct from singular values of other critical points.
\end{lemma}

\proof
It suffices to perturb each variable separately so each perturbation is a sum of 
monic polynomials of degree $d$ in a single variable $v_\ncount$. 
Since a polynomial $x_\ncount^d$ is versally unfolded by monomials of degree less than $d$,
we can choose each family of summands to consist of Morse functions only except for a finite number of
polynomials with a solitary degenerate critical point.

Since the critical points of the sum are the point such that each coordinate is a critical point of the corresponding
summand, we deduce, that to get a degenerate critical point at least in one coordinate the corresponding critical
point must be degenerate.

In fact it is non-generically degenerate if and only if either the critical point is degenerate in two coordinates
or its is non-generically degenerate in one coordinate.

But with our choice this happens only in a codimension two set.
Hence we can find a pencil as claimed in the $n$-dimensional family obtained from the one-dimensional
perturbation in each coordinate.
\qed

\begin{lemma}
\labell{a11gen}
There is a pencil of subdiagonally perturbed functions, such that at most two critical points have a common
critical value for all functions.
\end{lemma}

\proof
We consider the case $n=2$.
If in the first variable we pick a generic one-parameter family of polynomials,
there is a lower bound on the maximal distance of any three critical points.
In the second variable we pick a generic one-parameter family
with a polynomial in which all critical values are within this bound
of each other.

In the two-dimensional family, the line along this family does not meet the locus of multiple
conflicting values. So this locus must be contained either in codimension two, i.e.\ points, or in parallel lines.

Suppose there is such a parallel line and pick some point on it, then the conflicting critical points
must keep the same distance along that line, which is impossible since the corresponding 
discriminant curve is irreducible, so each pair can be moved to become arbitrarily close.

Hence there is no such parallel line, the locus of multiple conflicting values is in codimension two and thus
a suitable pencil of polynomials can be chosen.

In the general case $n\geq2$ we argue as before, but start with a generic pencil in all but the last variable.
\qed

The full genericity property we need comprises the two properties above and the property, that projection along
$L'$ is generic in the sense of lemma \reff{projv}.

\begin{lemma}
\labell{b-gen}
There is a generic line $L'$ in the base of the truncated subdiagonal unfolding.
\end{lemma}

\proof
By lemma \reff{projv} we know that there is a projection $p_v:\CC^{N-1}\to\CC^{N-2}$
such that a generic fibre is generic for $\bfami$. Notice that we have an open condition, hence a Zariski open
set of projections has that property.

Similarly the properties of the two preceding lemmas are open conditions.
Hence we conclude, that there is line $L'$ which meet simultaneously the conditions of the two lemmas and which
is a fibre of a projection with a property as in lemma \reff{projv}.
\qed

Now $L'$ can be used to compare the braid monodromy relations of \cite{habil}
with the relations imposed by geometric elements associated to $\bfami$ in
$\pi_1(\CC^{N}-\dfami)$.

To facilitate the use of citations from \cite{habil}
let us recall

\paragraph{Definition}
The \emph{braid monodromy group} is a subgroup of the automorphism group
of the free group generated by $t_i$. It is generated by the images of
all geometric elements associated to $\bfami$.

This subgroup is a subgroup of the braid group $\br$ considered as a group of automorphisms
by the Hurwitz action on the free group freely generated by elements $t_i$.

\begin{prop}
\labell{kernel}
Suppose $\pi_1(L'-\bfami_{L'})$ is generated by a geometric basis $r_b$ with base point close
to the Fermat pencil, and suppose $\Gamma$, the braid monodromy group of the Fermat polynomial,
is generated by $\{\b_s\}$,
then there is an identity of normally generated normal subgroups
$$
\langle\langle t_i\inv\phi_b(t_i)\rangle\rangle\quad=\quad\langle\langle t_i\inv\b_s(t_i)\rangle\rangle
$$
\end{prop}

\proof
From the previous results we infer, that the line $L'$ is generic also for the versal unfolding of the Fermat
polynomial. Hence the normal subgroup on the left hand side is the same as the
kernel of the map from the free group generated by the $t_i$ to the fundamental group of the discriminant
complement of a versal unfolding of the Fermat polynomial.

By the result of $\cite{habil}$ that kernel is given by the normal subgroup on the right hand side.
\qed

The explicit description of the normal subgroup in term of a finite set of relations will be taken from \cite{habil}
and used in the proof of the main theorem in section \reff{results}.

\section{Asymptotes}

Having fixed a preferred choice of generators for the discriminant complement in the previous section,
it is now the time to get the explicit relations imposed on our generators by the degenerations along the divisor $\afami$.

We first investigate the local situation to understand the impact of an asymptotic behaviour for
an adopted set of generators.

Our second aim is to relate the local situations which we spot at several transversal discs to $\afami$
via parallel transport to our distinguished set of generators.

Finally we have to establish that in such a way we have given sufficiently many relations in the sense
that every additional relation obtained from any local spot via any parallel transport is contained
in the normal closure of those we have listed.

\subsection{affine and local models}

\newcommand{\order}{n}
\newcommand{\Cpunc}{\CC-{\scriptstyle[}\order{\scriptstyle]}}  
\newcommand{\Cppunc}{\CC-[2]}  
\newcommand{\Dpunc}{D-{\scriptstyle[}m{\scriptstyle]}}  

To describe the geometry of the projected discriminant locally at a point of
$\afami$ we consider first the plane affine curves $C_\order$
\begin{eqnarray*}
C_\order\subset \CC^2 & : & y(yx^\order-1) \quad=\quad 0.
\end{eqnarray*}

We call $C_\order$ an \emph{asymptotic curve} of order $\order$,
since its equation can be rewritten as $y(y-\frac1{x^\order})$ to show that
the $y$-axis is a vertical asymptote and $\order$ is its pole order.

\begin{lemma}
\labell{asyp-I}
If in a vertical fibre $b_0,b$ are a geometric basis, then
the complement of a standard asymptotic curve $C_{\order}$ in $\CC^2$ has fundamental
group presentable as
\begin{eqnarray*}
\pi_1(\CC^2-C_{n}) & \cong & \langle b_0,a_0\,|\, b_0^n a_0 = a_0 b_0^n \rangle\\
& \cong &  \langle b_0,b\,|\, (bb_0)^n = (b_0b)^n \rangle,
\end{eqnarray*}
where $a_0$ links infinity.
\end{lemma}

\proof
The complement of $C_{n}$ is isomorphic to the complement of $C'_{n}$, 
a plane affine curve given by
$$
z(x^n-z),
$$
since both are complements of a quasi-projective curve in $\PP^1\times \CC$ given by
$$
yz(x^ny-z).
$$
The projection along $z$ exhibits the complement of $C'_{n}$ to be a fibre
bundle over the punctured complex line $\CC^*$
with fibre diffeomorphic to a $2$-punctured complex line $\Cppunc$.
The projection map has a section hence we get a short exact sequence with split surjection
$$
1\to \pi_1(\Cppunc)\to\pi_1(\CC^2-C'_{n})\to \pi_1(\CC-\{0\})\to 1.
$$
With base points $z_0>1,x_0=1$ we may choose elements
\begin{itemize}
\item
$b_0,a_0$ in $\pi_1(\Cppunc)$ represented by a geometric base of simple loops,
\item
$\bar c$ in $\pi_1(\CC-\{0\})$ with a lift $c\in \pi_1(\CC^2-C'_{n})$ via a section to $z=z_0$.
\end{itemize}
All relations follow from the conjugacy action of $c$ on the generators, which can be determined
as a braid monodromy and the triviality of $c$.

Indeed our claim corresponds to the conjugacy action given explicitly as
$$
c b_0 c\inv = (a_0b_0)^n b_0(a_0b_0)^{-n}.
$$
which is an immediate consequence of the braid monodromy given by $n$ full twists $\s_1^{2n}$.

The first presentation is then obtained by replacing $a_0$ by $(bb_0)\inv$ and rearranging.
\qed

Since $C_\order$ is invariant under $(x,y)\mapsto (\l^\order y,\l\inv x)$
its complement can be strongly retracted to any fibred neighbourhood of
the $y$-axis.

\paragraph{Definition}
A curve $C\subset D\times\CC$ is called an \emph{asymptotic germ} if
\begin{enumerate}
\item
the fibration to $D$ is locally trivial onto the punctured disc $D^*$,
\item
the projective closure of $C$ in $\PP^1\times D$ is regular.
\end{enumerate}

\paragraph{Example}
Consider a plan affine curve defined by
$$
(y^m-1)(x^ny-1).
$$
The restriction to a sufficiently small neighbourhood of $x=0$
yields an asymptotic germ $C_{n,m}$.

\begin{lemma}
\labell{asyp-II}
With the embeddings of $D\times A$ as a tubular collar of the component $y=0$ into
$(D\times\CC)-C_{n}$ and as an exterior collar of $(\Dpunc)\times \CC$:
$$
(D\times\CC)-C_{n,m} = (D\times\CC-C_{n})\cup_{D\times A}(\Dpunc)\times\CC.
$$
\end{lemma}

\begin{prop}
\labell{local1}
If $C$ is an asymptotic curve of fibre degree $m+1$ and pole order $n$ then
there is a geometric basis $b_0,b_1,...,b_m$ in a fibre of the complement such that
$$
\pi_1(\CC\times D-C) \quad \cong \quad \langle b_0,b_1,...,b_m\,|\, b^n b_0=b_0 b^n\rangle
$$
where $b=b_m\cdots b_1$.
\end{prop}

\proof
The curve $C$ is fibre-isotopic in $\CC\times D$ to the model $C_{n,m}$, so we can employ
the decomposition given above.
We choose a geometric basis such that $b_0$ is supported in the outside part and
the $b_1,...,b_m$ in the inside disc.
Then $\pi_1((\Dpunc)\times\CC) $ is the free group generated freely by $b_1,...,b_m$
and $\pi_1(D\times \CC- C_n)$ is isomorphic to $\pi_1(\CC^2-C_n)$ and hence presented as
$$
\langle b_0,a \,|\,(b_0a)^n = (a b_0)^n \rangle,
$$
where $a$ is the geometric generator around infinity.

We can the apply the vanKampen theorem with the observation that $b=b_m\cdots b_1$ and
$b_0a$ represent the same generator of the fundamental group in the intersection, so
\begin{eqnarray*}
\pi_1 & \cong &
\langle b_0,b_1,...,b_m,a \,|\, a\inv = bb_0,\, (b_0 a)^n=(a b_0)^n \rangle\\
& \cong &
\langle b_0,b_1,...,b_m \,|\, (b_0 b_0\inv b\inv)^n = (b_0\inv b\inv b_0)^n \rangle\\
& \cong &
\langle b_0,b_1,...,b_m \,|\, b^{-n} = b_0\inv b^{-n}b_0 \rangle
\end{eqnarray*}
from which the claim is immediate.
\qed

A geometric basis as in the proposition is called \emph{adopted} to the asymptote.
and provides a tool to make the second set of relations in lemma \ref{pd}
more explicit.

\begin{lemma}
\labell{asyp}
Suppose a set of generators $r_a$ for $\pi_1(\CC^{N-1}-\afami)$ is given
such that each $r_a$ is associated to an asymptotic degeneration of order $n$,
then there are two sets of relations describing the same normally generated subgroup:
$$
\langle\langle \rho_a t_i\inv\rho_a\inv\phi_a(t_i)\rangle\rangle
\quad=\quad
\langle\langle (\tau_a\inv\d_0)^n(\d_0\tau_a\inv)^{-n}\rangle\rangle)
$$
where $\tau_a$ is the simple geometric element corresponding to $b_0$ in the local model.
\end{lemma}

\proof
The claim is immediate from the fact that the family of a transversal disc to $\afami$ is topological
equivalent to the local model constructed above, so both ways to construct relations have to
yield the same normal subgroup.
\qed

\subsection{a core subfamily}

In this paragraph we consider a suitable family $\gfami$ of polynomials in which our parallel transport
will take place.
\begin{eqnarray*}
& & \z\bigg( x_n^d-d\l_n^{d-1}x_n+
\sum_{\ncount=1}^{n-1}\left(x_\ncount^d-d\l_\ncount^{d-1}x_\ncount\right)\bigg)\\
& + & \eta\bigg((d-1)^2x_n^d-d(d-1)\l_nx_n^{d-1} +
\sum_{\ncount=1}^{n-1}\left((d-1)x_\ncount^d-d\l_\ncount x_\ncount^{d-1}\right)\bigg)\\
& - & \l d(d-1)\sum_{\ncount=1}^{n-1}\l_\ncount^{d-1} x_\ncount x_n^{d-1}
\end{eqnarray*}
parameterised by $\z,\eta,\l,\l_1,...,\l_n$. We study some of its properties which will serve later to choose
appropriate paths of polynomials.

To highlight te distinguished role given to $x_n$ we replace this variable by $y$ in this instance.

\begin{lemma}
\labell{crit1}
The point $(\l_1,...,\l_n)$ is a critical point of each of the polynomials in the restricted
family $\gfami|_{\l=0}$.
\end{lemma}

\proof
We need the components of the gradient with respect to the variables $x_\ncount$, where $y=x_n$
has to be handled separately.
\begin{eqnarray*}
& \del_\ncount: & \z\left( dx_\ncount^{d-1}-d\l_\ncount^{d-1}\right)
+ \eta\left(d(d-1)x_\ncount^{d-1}-d(d-1)\l_\ncount x_\ncount^{d-2}\right)\\
& \del_y: & \z\left( dy^{d-1}-d\l_n^{d-1}\right)
+ \eta\left(d(d-1)^2y^{d-1}-d(d-1)^2\l_n y^{d-2}\right)
\end{eqnarray*}
Then it is easy to see, that the gradient vanishes at the given point.
\qed

\begin{lemma}
\labell{maxmod1}
In the restricted family $\gfami|_{\l=0}$ with positive real parameters $\z,\eta$
the critical point $(\l_1,...,\l_n)$ is the critical points with critical value of largest modulus.
\end{lemma}

\proof
Our first claim is that any coordinate $x_\ncount$ of a critical point either equals $\l_\ncount$
or has modulus distinct from $|\l_\ncount|$.

So suppose $x_\ncount$ is not $\l_\ncount$ but of the same modulus, that is
$x_\ncount=\rho\l_\ncount$ for some complex number of modulus one distinct from 1.
Then the vanishing of the gradient implies
\begin{eqnarray*}
\z(\rho^{d-1}-1)\l_\ncount^{d-1}+\eta(d-1)(\rho^{d-1}-\rho^{d-2})\l_\ncount^{d-1} & = & 0\\
\implies\qquad
\z(\rho^{d-1}-1)+\eta(d-1)(\rho^{d-1}-\rho^{d-2}) & = & 0.
\end{eqnarray*}
But the two differences of complex numbers of modulus one are not parallel. In fact
the last equation is solvable only for $\rho=1$ or $\eta=0$.
\\

Next we deduce that any coordinate $x_\ncount$ of a critical point has modulus less or equal 
to $|\l_\ncount|$. That follows from the connectedness of the space of parameters, and the fact that
the claim is true for $\eta=1$ and $\z$ sufficiently small.

No coordinate may then pass the threshold of modulus $|\l_\ncount|$, since otherwise we would
get a contradiction to our first claim.
\\

Finally we can make an estimate for the critical values.
Imposing the gradient condition the function coincides on the critical points with
$$
\z(1-d)\left(\l_n^{d-1}y+\sum_{\ncount=1}^{n-1}\l_\ncount^{d-1}x_\ncount\right)
-
\eta\left((d-1)\l_ny^{d-1}+\sum_{\ncount=1}^{n-1}\l_\ncount x_\ncount^{d-1}\right)
$$
Hence a bound is given by
$$
\z(d-1)\left(\l_n^{d}+\sum_{\ncount=1}^{n-1}\l_\ncount^{d}\right)
+
\eta\left((d-1)\l_n^{d}+\sum_{\ncount=1}^{n-1}\l_\ncount^{d}\right)
$$
which is strict in all cases except for the critical point $(\l_1,...,\l_n)$, when
it is attained.
\qed

\begin{lemma}
\labell{crit2}
There is some $Z$ with $Z=1-\l^d(\sum \l_\ncount^d)^{d-1}$ up to terms of order at least one in $\l_n$
such that $Z\inv(\l_1,...,\l_{n-1},\l_nZ+\l\sum\l_\ncount^d)$ is a critical point
in the restricted family $\gfami_{\z=0,\eta=1}$.

Moreover $Z$ is real with $0<Z<1$ if the $\l_\ncount^{d}$ are real and sufficiently small.
\end{lemma}

\proof
We need the gradient components of the restricted family:
\begin{eqnarray*}
& \del_\ncount: &
d(d-1)x_\ncount^{d-1}-d(d-1)\l_\ncount x_\ncount^{d-2}-\l d(d-1)\l_\ncount^{d-1}y^{d-1}\\
& \del_y: &
d(d-1)^2y^{d-1}-d(d-1)^2\l_ny^{d-2}-\l d(d-1)^2y^{d-2}\sum\l_\ncount^{d-1}x_\ncount
\end{eqnarray*}
The second line vanishes for the given point, whatever $Z$ is.
For each $\ncount$ the first line gives an expression which is a multiple of $\l_n$ if
we put $Z=1-\l^d(\sum \l_\ncount^d)^{d-1}$.
\\

The claim $Z$ real and positive is obvious. To prove $Z<1$ we note that $Z\geq1$ in 
any gradient component of the first kind leads to a contradiction.
\qed

\begin{lemma}
\labell{maxmod2}
Suppose that the parameters $\l_\ncount$ are fixed such that each $\l_\ncount^d$ is a positive real number.
Then in the restricted family $\gfami|_{\z=0}$ with sufficiently small positive real parameters $\l,\l_n$
the critical point $(\l_1,...,\l_n)$ of largest modulus varies continuously.
\end{lemma}

\proof
For $\l=0$ the claim follows from lemma \reff{maxmod1}. For $\l=\l_n=0$ the distinguished critical
point has higher multiplicity and has $y=0$.

Since other critical values are still of smaller modulus, we have to check only those critical points
which limit to the distinguished one for $\l,\l_n\to 0$.

From the $\del_y$ part of the gradient we deduce that any such critical point distinct
from the distinguished one has constantly $y=0$ and accordingly $x_\ncount=\l_\ncount$.

A brief check against the distinguished critical value with the properties as given by lemma \reff{crit2}
proves our claim.
\qed

\begin{lemma}
\labell{maxmod3}
Suppose that the parameters $\l_\ncount$ are fixed such that each $\l_\ncount^d$ is a positive real number.
Then for positive $Z=1-\l^d(\sum \l_\ncount^d)^{d-1}$
the restricted family $\gfami|_{\z=0,\l_n=0}$ with positive real parameter $\l$
has the critical point $Z\inv(\l_1,...,\l_{n-1},\l\sum\l_\ncount^d)$ with critical value
of largest modulus.
\end{lemma}

\proof
For any critical point $y=0$ or $y=\l\sum\l_\ncount^{d-1}x_\ncount$.
In the first case it is immediate that all occurring $x_\ncount$ are bounded by that of the
distinguished critical point. Hence also the critical value is bounded.

In the second case let us write each coordinate as
$$
x_\ncount=\rho_\ncount\frac{\l_\ncount}{Z}.
$$
Suppose $\rho_1$ is of largest modulus among the $\rho_\ncount$. We derive 
from the corresponding gradient component:
$$
\frac{Z}{\rho_1}\quad=\quad
1-\l^d\left(\sum\l_\ncount^d\frac{\rho_\ncount^d}{\rho_1^d}\right)^{d-1}.
$$
By hypothesis the right hand side has modulus at least $Z$ with
equality only if all $\rho_\ncount$ are equal. Hence on the left hand side
the modulus of $\rho_1$ may be at most $1$.

So in fact our claim follows.
\qed

\subsection{asymptotic arcs, paths, and induced paths}

We recall that the space $\CC^{N-1}$ parameterises hypersurface pencils containing
the multiple hyperplane $x_0^d$ or -- equivalently -- polynomials of degree at most $d$
in $\CC[x_1,...,x_n]$ with vanishing constant coefficient.

Hence arcs and paths in $\CC^{N-1}$ are given by families of polynomials parameterised by a real
interval.
\\

We introduce the following arcs starting at the base point
$$
\z=1,\,\eta=0,\,\l=0,\,v_\ncount.
$$
\begin{enumerate}
\item
A set of arcs in bijection to tuples $\jbold=j_1,...,j_{n}\in \ind$ where $\l_\ncount$ is chosen such
that $\l_\ncount^{d-1}$ moves constantly along a circle segment
from $v_\ncount$ to $v_\ncount\xi^{1-j_\ncount}$, cf. the construction of the paths $\oo_\jbold$
in \reff{pathsandhomotopies}.

The endpoints we denote by $v(\jbold)$.

\item
A set of arcs with all parameters fixed, such that $\l=0$, $\l_\ncount^{d-1}=v_\ncount(\jbold)$
except for the convex relation $\z+\eta=1$,

\item
A set of arcs with all parameters fixed at some endpoint of the previous set,
except for $\l_n$ and $\l$. Along the arc $\l_n$ vanishes and $\l$ increases to some
sufficiently small positive real value.

\item
A set of arcs with all parameters fixed at some endpoint of the previous set,
except for $\l$ which increases to $\frac{1}{\sum\l_\ncount^d}$.
\end{enumerate}

For each tuple
$\jbold=j_1,...,j_{n}$ we get exactly one composite connected arc.
Then the composite arc $\a(\jbold)$ leads from the Fermat point to a point on $\afami$ without intersecting
$\afami$ elsewhere.

Hence there exists a corresponding path $\wp(\jbold)$ closed at the Fermat point which is
arbitrarily close to $\a(\jbold)$ in the complement of $\bfami$ and $\afami$ and a geometric element for
$\afami$.
\\

\begin{prop}
\labell{transport}
The relation imposed on the generators along the path $\wp(\jbold)$
is 
$$
(t_{\jbold}\inv\d_0)^{d-1}=(\d_0t_{\jbold}\inv)^{d-1}.
$$
\end{prop}

\proof
In fact the relation is given by the local relation from prop. \reff{local1}. The point
is that we have to identify the local element in terms of the chosen generators
at the base point.

Now the difficult construction of this section shows, that parallel transport first moves a
representative of $t_\jbold$ to the representative $t_{\ibold_0(1)}(\jbold)$. And then this representative is transported to
$b_0$, since the inner puncture remains always the puncture of largest modulus.

Of course the loop around all punctures is transported to the loop around all punctures,
hence parallel transport identifies $\d_0$ with $bb_0$.
We conclude $b=\d_0t_\jbold\inv$ and get
\begin{eqnarray*}
(\d_0t_\jbold\inv)^{d-1}t_\jbold & = & t_\jbold(\d_0t_\jbold\inv)^{d-1}\\
\iff
t_\jbold\inv\d_0( t_\jbold\inv\d_0)^{d-2}& = & (\d_0t_\jbold\inv)^{d-1},
\end{eqnarray*}
which yields the claim.
\qed

\begin{prop}
\labell{generators1}
The set of paths $\wp(\jbold), \jbold=\jbold'1, \jbold'\in I_{n-1,d}$ generates $\pi_1(\CC^{N-1}-\afami)$.
\end{prop}

\proof
It suffices to check the claim for the projection corresponding to $x_0=0$ and $x_n=1$.
Then only the parts of $\wp(\jbold)$ corresponding to the first and the last arc component
are non-trivial. 

We find that the parts of $\wp(\jbold)$ at the base point coincide with the path $\oo(\jbold)(t_{\ibold_0(1)})$
only in the middle part the coefficients of $x_\ncount x_n^{d-1}$ are increased, instead of
decreasing the coefficient of $x_n^{d}$.
But this gap is easily overcome by a smooth homotopy which rescales $x_n$.

Hence the projected paths are homotopic to the set of $\oo(\jbold)(t_{\ibold_0(1)})$
and therefore to the set $t_\jbold, \jbold=\jbold'1, \jbold'\in I_{n-1,d}$
which is known to generate the fundamental group, since $\CC^{N-1}-\afami\simeq
\CC^{N'}-\dfami_{n-1,d}$.
\qed

\section{projective quotient}
\label{projective}

It is time now to recall that our interest is in the fundamental group of the
discriminant complement $\und$.
The relation to our previous results can be expressed as follows.

\begin{lemma}
\labell{projquot}
The projective discriminant complement $\ufami$ is the quotient of
the complement $\CC^N-\dfami$ to the affine discriminant cone by the
diagonal $\CC^*$-action
and for all $u\in\CC^N-\dfami$
there is an exact sequence
$$
1\tto\pi_1(\CC^*,1)\tto\pi_1(\CC^N-\dfami, u),\tto\pi_1(\ufami,[u])\tto 1.
$$
\end{lemma}

\proof
Since the discriminant is an affine cone which is preserved by the $\CC^*$ action
so is the complement.
The quotient map then naturally gives rise to a fibration $\CC^N-\dfami$ to $\und$ with
fibre $\CC^*$.

Its long exact homotopy sequence provides an exact sequence
$$
\pi_1(\CC^*,1)\tto\pi_1(\CC^N-\dfami, u),\tto\pi_1(\ufami,u)\tto \pi_0(\CC^*).
$$
Of course  $\pi_0(\CC^*)$ is the trivial group, so to get the sequence of our claim
it remains to prove that the first map is injective.

For that we look at the abelianisation
of $\pi_1(\CC^N-\dfami, u)$ which is isomorphic to $\ZZ$ because of the irreducibility of
the discriminant.
In fact the abelianisation map is given by the linking number.

Since the image of a non-trivial closed path in $\CC^*$ links the discriminant
non-trivially, the claim is immediate.
\qed

To proceed we pick $u$ to be the Fermat point corresponding to the Fermat hypersurface 
$$
u:\quad x_0^d+x_1^d+\cdots+x_n^d.
$$
Its $\CC^*$ orbit belong entirely to the Fermat family $\ffami$ to which we shift our attention for the moment
$$
a_0x_0^d+a_1x_1^d+\cdots+a_nx_n^d.
$$

\paragraph{Definition}
The element $\delta_\ncount$ is defined as the element represented by the path in the Fermat family
$\ffami$ based at $(1,...,1)$ which is constant in all components except for $a_\ncount=e^{it}$.

\begin{lemma}
\labell{deltaab}
The $\delta_\ncount$ all commute with each other and
$$
\prod\d_\ncount\quad=\quad 1\in \pi_1(\und).
$$
\end{lemma}

\proof
The Fermat family $\ffami$ has discriminant given by the normal crossing
divisor given by the coordinate hyperplanes and hence the fundamental group of the complement is
abelian.

Since the $\d_\ncount$ are geometric generators associated to the $n+1$ hyperplanes
the fundamental group of the projectivised complement is the free abelian group
generated by the $\d_\ncount$ modulo
the subgroup generated by their product.

Of course this relation maps homomorphically to $\pi_1(\und)$.
\qed

\paragraph{Remark}
We are going to combine results on fundamental groups of $\und$ and various of
its subspaces which do \emph{not} have the same base point.

But in fact all these base points are contained in a ball in $\und$ so our convention is
that all occurring fundamental groups are identified using a connecting path
for their base points inside this ball, which makes the identification 
unambiguous.
\\

To formulate our results we employ enumeration functions
\begin{eqnarray*}
\ibold_\ncount &: & \{1,...,(d-1)^n\}\to \ind
\end{eqnarray*}
which are most conveniently described by the total order $\kord$ they induce on $\ind$.

The enumeration function $i_0$ has been defined above as the standard lexicographical order of 
$\ind$ derived from the \emph{reverse} order of each factor $\{1,...,d-1\}$ of $\ind=\{1,...,d-1\}^n$.

The modified order $\kord$ of $\ind$ is again a lexicographical order, this time
each factor is ordered reversely -- that is by $>$ --
except for the $\ncount$-th factor which is naturally ordered, from small to large.

For $n=1$ in particular $\prec_1$ is the natural order $<$.
\\

By lemma \reff{cox} the element $\d_0$ can be expressed in the
geometric basis $t_\ibold$ of $\pi_1(L_0-\dfami_L)$
using the enumeration function $\ibold_0:\{1,...,(d-1)^n\}\to \ind$:
$$
\d_0\quad =\quad\prod_{m=1}^{(d-1)^n} t_{\ibold_0(m)}.
$$

From this we proceed to get expressions for the $\d_\ncount$ as well:

\begin{lemma}
\labell{anfang}
In case $n=1$ the element $\d_1$ can be expressed in the geometric basis $t_\ibold$
using the enumeration $i_{1}:\{1,...,d-1\}\to I_{1,d}$:
$$
\d_1\quad=\quad \prod_{m=1}^{d-1} t_{i_{1}(m)}\quad =\quad t_1t_2\cdots t_{d-1}.
$$
\end{lemma}

\proof
This is the case studied by Zariski. We only need to combine our result $\d_0\d_1=1$, lemma
\reff{deltaab}, with $\d_0\,=\, t_{d-1}\cdots t_1$ and 
his result
$$
t_1\cdots t_{d-1}t_{d-1}\cdots t_1\quad=\quad 1.
$$
to conclude $\d_1\,=\, t_1\cdots t_{d-1}$
as claimed.
\qed

\begin{lemma}
\labell{dim-ind}
In case of general $n$ the element $\d_1$ can be expressed in the geometric basis $t_\ibold$
using the enumeration $\ibold_{1}:\{1,...,(d-1)^n\}\to \ind$ as
$$
\d_1\quad=\quad \prod_{m=1}^{(d-1)^n} t_{\ibold_{1}(m)}.
$$
\end{lemma}

\proof
We consider the space $\ufami_{1,d}$ and a homotopy $H$ between the representing path of 
$\d_1$ and the representing path of $t_1\cdots t_{d-1}$, cf. lemma \reff{anfang}.

We consider the map defined on $\ufami_{1,d}$ to $\psym$ defined by
$$
f\quad \mapsto\quad f+(x_2^d+x_3^d+\cdots+x_n^d)
$$
which maps $H$ to $\und$.

There $H$ provides a homotopy between $\d_1$ and a path representing the product of the images of
$t_1,t_2,...,t_{d-1}$. 
These images are exactly the factors of $\d_0$ which are formed by all generators with the same first index
component.
Hence the claimed relation holds.
\qed

\begin{lemma}
\labell{ind-anf}
Consider the Hefez Lazzeri family
$$
x_1^d-v_1x_1x_0^{d-1}+x_2^d-v_2x^2x_0^{d-1}+zx_0^d.
$$
Suppose $t_\ibold$ and $t_\ibold'$ are geometric generators for real parameters $v_1\ll v_2$ respectively
$v_1\gg v_2$ of sufficiently distinct magnitude,
then there is a path connecting the respective base points such that the associated isomorphism
on fundamental groups is given by
$$
t_{i_1i_2}\quad \mapsto \quad t_{i_2i_1}'\quad 
$$
\end{lemma}

\proof
We first convince ourselves that $t_{11}=t'_{11}$ which follows immediatly if we change $v_1,v_2$ continuously
in the real line swapping places since the extremal real puncture will keep that property and hence the
corresponding geometric element will not be changed.

To move $t_{i_1i_2}$ we first proceed along a path $\oo$ so that it becomes the $t_{11}$ in the new system,
then do the same as above and finally employ a path $\oo$ again to come to the final position.
\qed

The argument of the proof readily generalises to the case of more variables.

\begin{lemma}
\labell{ind-ind}
Consider the Hefez Lazzeri family
$$
\sum (x_\ncount^d-v_\ncount x_\ncount x_0^{d-1})+zx_0^d.
$$
Suppose $\pi$ is a permutation such that $t_\ibold$ and $t_\ibold'$ are geometric generators
 for real parameters $v_1\ll v_2\ll\cdots\ll v_n$ respectively
$v_{\pi(1)}\ll v_{\pi(2)}\ll \cdots\ll v_{\pi(n)}$ of sufficiently distinct magnitude
then there is a path connecting the respective base points such that the associated isomorphism
on fundamental groups is given by
$$
t_{i_1i_2\cdots i_n}\quad \mapsto \quad t_{i_{\pi(1)}i_{\pi(2)}\cdots i_{\pi(n)}}'
$$
\end{lemma}

\begin{lemma}
\labell{proj-I}
In case of general $n$ an expression for $\d_\ncount$ is given by
$$
\prod_{m=1}^{(d-1)^n}t_{\ibold_\ncount(m)}\quad=\quad 1.
$$
\end{lemma}

\proof
We get the expression for general $\d_\ncount$ using a transposition $\pi=(1\ncount)$ in the previous
lemma on the expression for $\d_1$.

Then it is obvious that the non-reversed order is transfered from the first entry into the $\ncount$-th entry.
\qed

\begin{prop}
\labell{projrel}
The image of a generator of $\pi_1(\CC^*)$ in the fundamental group $\pi_1(\und)$
for the natural map to the orbit of the base point is given in the Hefez-Lazzeri geometric generators
$t_\ibold$ as
$$
\prod_{\ncount=0}^n\:\prod_{m=1}^{(d-1)^n}t_{\ibold_\ncount(m)}\quad=\quad 1.
$$
\end{prop}

\proof
Immediate from the previous.
\qed

It should be remarked that our formulation of results does not respect the symmetry of the geometry, but this 
was to be expected as we have broken the symmetry from the begining by picking $x_0$ to be special.

\section{Conclusion}
\label{results}

Finally we are ready to prove the Main Theorem. As promised before we present the necessary result from \cite{habil}
first.
There we express the fundamental group of the complement to the discriminant in any versal unfolding of any
Brieskorn-Pham singularity.
In the special case of Fermat polynomials $x_1^d+\cdots x_n^d$ the result can be most conveniently
presented in term of an integer lattice with basis $\{\vect_i,i\in I_{n,d}\}$ and bilinear form defined by
$$
\langle \vect_i,\vect_j\rangle \quad=\quad\left\{
\begin{array}{cl}
0 & \text{if } |i_\nu-j_\nu|\geq2\text{ for some }\nu,\\
0 & \text{if } (i_\nu-j_\nu)(i_\mu-j_\mu)<0 \text{ for some } \nu,\mu,\\
-2 & \text{if } i=j\\
-1 &  \text{ else}.
\end{array}
\right.
$$

\begin{thm}[\cite{habil}]
For the complement of the discriminant of a versal unfolding of the singularity $x_1^d+\cdots x_n^d$:
$$
\pi_1\quad\cong\quad
\left\langle t_i,\,i\in I\,\left|
\begin{array}{ccl}
t_it_j=t_jt_i & \text{if} & \langle \vect_i,\vect_j\rangle=0\\
t_it_jt_i=t_jt_it_j & \text{if} & \langle \vect_i,\vect_j\rangle\neq0\\
t_it_jt_kt_i=t_jt_kt_it_j & \text{if} & \langle \vect_i,\vect_j\rangle
\langle \vect_j,\vect_k\rangle\langle \vect_k,\vect_i\rangle\neq0
\end{array}
\right.\right\rangle.
$$
\end{thm}

The lattice itself is naturally encoded in the intersection graph $\Gamma_{n,d}$ on the
set of vertices $V(\Gamma)=\ind$ with set of edges
$$
E(\Gamma)\quad=\quad \{(i,j)\,|\,\langle\vect_i,\vect_j\rangle\neq0\}.
$$
Example of these graphs we have given in the introduction. 

The Main Theorem is thus
made precise with the definition of $\Gamma_{n,d}$ just given and the enumeration
function $\ibold_\ncount$ of the previous section defining the distinguished elements
$\d_\ncount$
With the notation developed throughout the paper we can give a changed and more
precise statement of our main result:

\paragraph{Main Theorem}
\unitlength=1pt
The complement of the discriminant
for degree $d$ hypersurfaces in $\PP^n$ has $\pi_1$ given by quotient of
$$
\left\langle t_i,\,i\in I\,\left|
\begin{array}{cll}
 t_i t_j=t_j t_i,& \text{for }
 {}_i\!\cdot \phantom{\text{---}} \cdot_j\,\\
 t_i t_j t_i=t_j t_i t_j,&  \text{for } {}_i\!\cdot\, \text{---} \cdot_j,\\
 t_it_jt_kt_i=t_jt_kt_it_j, & \text{for }
{}_i\,
\begin{picture}(20,5)
\put(0,4){\circle*{1}}
\put(20,4){\circle*{1}}
\put(10,-6){\circle*{1}}
\put(5,-9){${}_j$}
\bezier{20}(4,4)(10,4)(16,4)
\bezier{20}(16,4)(14,4)(12,4)
\bezier{30}(2,2)(5,-1)(8,-4)
\bezier{30}(18,2)(15,-1)(12,-4)
\end{picture}
\,{}_k\\[2mm]
\end{array}
\right.\right\rangle.
$$
by the normal subgroup generated by the relations
\begin{eqnarray}
(t_i\inv\delta_0)^{d-1} & = & (\delta_0t_i\inv)^{d-1}\\
\delta_0\delta_1\cdots\delta_n & = & 1
\end{eqnarray}

\proof
We first show, that the fundamental group of the complement to the affine cone of the
discriminant is presented as in the claim except for the last relation (5).

By lemma \ref{b-gen} there is a fibre $L'$ for some projection $p_v:\CC^{N-1}\to\CC^{N-2}$
such that lemma \ref{genab} applies;
elements $r_b$ of a geometric basis for $L'$ and elements $r_a$ generating $\pi_1(\CC^{N-1}-\afami)$
form a generating set for $\pi_1(\CC^{N-1}-\afami-\bfami)$.

So in principle a presentation can be obtained with the help of lemma \ref{pd}.

By prop.\ \ref{generators1} we may choose to take the generators $r_a$ to be represented by the
paths $\wp(\jbold)$. Then the relations of second type in \ref{pd} can be replaced
using those of prop.\ \ref{transport}, so by lemma \ref{asyp} we get the relations in (4).

From prop.\ \ref{kernel} we infer that we may replace the first set of relations by any other
set normally generating the kernel of the map from the free group
to the fundamental group of the discriminant complement in the singularity unfolding,
in particular by the set given in the above result from \cite{habil}.

Proceeding now to the projective quotient it suffices by lemma \ref{projquot} to
add just one relation, in fact that in (5) as we have shown in prop.\ \ref{projrel}.
\qed

In fact prop.\ \ref{generators1} allows us to use less relations, but we have kept them for symmetry.

\subsection{Remarks}

In cases where the fundamental groups have been described before,
our presentation is different.

\paragraph{Example}
\labell{slz}
In the case of cubic curves we get a group generated by four elements $t_1,t_2,t_3,t_4$
subject to the relations
\begin{itemize}
\item
$t_2t_3=t_3t_2$,
\item
$t_it_jt_i=t_jt_it_j$ if $(ij)\in\{(12),(13),(24),(34)\}$,
\item
$t_it_jt_kt_i=t_jt_kt_it_j$ if $(ijk)\in\{(124),(134)\}$,
\item
$t_4t_3t_2t_4t_3t_2t_1=t_1t_4t_3t_2t_4t_3t_2$,
\item
$t_3t_2t_1t_3t_2t_1t_4=t_4t_3t_2t_1t_3t_2t_1$,
\item
$t_4t_3t_2t_1t_2t_1t_4t_3t_3t_1t_4t_2=1$
\end{itemize}
The same group is described by Dolgachev and Libgober \cite{dl} as an extension of $\slz$ by a
$27$ element Heisenberg group over $\ZZ/3\ZZ$. It may be interesting to give an isomorphism
to the presentation above, and to spot the torsion element in particular.

\paragraph{Example}
In the case of cubic surfaces we get a presentation by $8$ generators $t_1,...,t_8$ subject
to the relations
\begin{itemize}
\item
$t_it_j=t_jt_i$ if $(ij)\in\{(23),(25),(27),(35),(36),(45),(46),(47),(67)\}$,
\item
$t_1t_jt_1=t_jt_1t_j$ for $j\in\{2,3,...,8\}$,
\item
$t_it_jt_kt_i=t_jt_kt_it_j$ for $(ijk)$ in
\begin{quote}
$
\{(124),(126),(128),(134),(137),(138),(148),(156),(157),(158),$\\
$(168),(178),(248),(268),(348),(378),(568),(578)\}
$
\end{quote}
\item
$t_i(t_i\inv\d_0)^{2}=(t_i\inv\d_0)^{2}t_i$, for $i\in\{1,2,...,8\}$ and $\d_0=t_8\cdots t_1$,
\item 
$t_8t_7t_6t_5t_4t_3t_2t_1 t_4t_3t_2t_1t_8t_7t_6t_5 t_6t_5t_2t_1t_8t_7t_4t_3 t_7t_5t_3t_1t_8t_6t_4t_2=1$
\end{itemize}
In fact we have neglected some of the redundant relations from the list in the Main Theorem.

Our presentation seems to be akin to that given by Looijenga \cite{lo} and again it would be nice
to have an explicit isomorphism.

\subsection{$\pi_1$ of complements to partially smooth discriminant curves}

We finally turn to the proof of proposition \reff{tribute}. For the convenience of the reader we recall:

\paragraph{Definition}
A reduced plane curve $C_0$ is called a limiting curve of a curve $C_1$ if there is a
smooth family $C_t$
of plane curves which is equisingular for $t\neq0$.

$C_1$ is called a partial smoothing of $C_0$ if they are not equisingular.
\\

\proof[ of prop.\ \reff{tribute}]
In case $C_0$ is a generic plane section of any discriminant $\dfami_{n,d}$ its only singularities
are ordinary cusps and nodes. Hence a partial smoothing $C_1$ must contain one node or one
cusp less than $C_0$.
Hence either a braid relation or a commutation relation must be replaced by an identity relation.

The essential observation is, that in all cases the cuspidal edge
and the Maxwell stratum are irreducible sets.
That implies that replacing any braid relation or commutation relation by an identity relation
must result in replacing all braid or all commutation relations by the corresponding identity relation.

Hence we can do this on the particular presentations we have for the complements of $\dfami$ and $C_0$.
Since each graph $\Gamma$ is connected, in case of a cusp smoothing all generators are identified.
Hence the group is cyclic and its degree is given by the length of the relation (5), since all other relations
become trivial.

In the case of a node smoothing we try to get the same result.
As in the argument of Zariski \cite{z} we need a set of three generators which only support a single edge of the graph.
In the Zariski case $n=1$ this is provided by generators $t_1,t_2,t_4$ for $d\geq5$,
in case $n=2$ by $t_{11},t_{12},t_{31}$ for $d\geq 4$ and
in case $n\geq 3$ by $t_{112...},t_{121...},t_{212...}$ for $d\geq3$.

The two commutation relations force all three elements to be identified in the presence of a node
smoothing and therefore the third relation, a braid relation, is also replaced by an identification
of generators, hence we may conclude as in the case of a cusp smoothing.
\\

Since for $d=1,2$ the curve $C_0$ is smooth the only cases where smoothing of a node does not force the
residual fundamental group to be smooth are the special cases of $(n,d)=(1,3)$ and $(n,d)=(2,3)$.

The first of these has been handled by Zariski already.
In the case $(2,3)$ we get a duodectic which is a smoothing of a generic plane section of $\dfami_{2,3}$.
From the presentation of the corresponding fundamental group above we get a presentation of
$\slz$ as the quotient of the braid group $\br_3$ by the square of the central element, just by imposing
$t_3=t_2$ and $t_4=t_2t_1t_2\inv$:
$$
\langle t_1,t_2\;|\; t_1t_2t_1=t_2t_1t_2, (t_2t_1)^6=1\,\rangle,
$$
and other relations becoming trivial.
\qed


\end{document}